\documentclass[a4paper,11pt,onecolumn,twoside]{article}
\usepackage{mathrsfs}
\usepackage{mathrsfs}
\usepackage{amsfonts}
\usepackage{fancyhdr}
\usepackage{amsmath,amsfonts,amssymb,graphicx,amsthm}
\usepackage{cite}
\allowdisplaybreaks
\begin{document}

\title{\bf Global Existence and Vanishing Dispersion Limit of Strong/Classical Solutions to  the One-dimensional  Compressible  Quantum Navier-Stokes Equations with Large Initial Data}
\author{{\bf Zhengzheng Chen}\\
Center for Pure Mathematics, School of Mathematical Sciences, Anhui University,\\ Hefei 230601, China\\[2mm]
{\bf Huijiang Zhao}\thanks{Corresponding author. E-mail:
hhjjzhao@hotmail.com}\\
School of Mathematics and Statistics, Wuhan University, Wuhan 430072, China\\
and\\
Computational Science Hubei Key Laboratory, Wuhan University, Wuhan 430072, China}

\date{}

\vskip 0.2cm

\maketitle

\vskip 0.2cm \arraycolsep1.5pt
\newtheorem{Lemma}{Lemma}[section]
\newtheorem{Theorem}{Theorem}[section]
\newtheorem{Definition}{Definition}[section]
\newtheorem{Proposition}{Proposition}[section]
\newtheorem{Remark}{Remark}[section]
\newtheorem{Corollary}{Corollary}[section]

\begin{abstract}
We are concerned with  the global existence and vanishing dispersion limit of strong/classical solutions to the Cauchy problem of the  one-dimensional  barotropic compressible quantum  Navier-Stokes equations, which  consists of the compressible Navier-Stokes equations with a linearly density-dependent viscosity and a nonlinear third-order differential operator known as the quantum Bohm potential. The pressure $p(\rho)=\rho^\gamma$ is considered with $\gamma\geq1$ being a constant. We focus on the case when the viscosity constant $\nu$ and the Planck constant $\varepsilon$ are not equal. Under some suitable assumptions on $\nu,\varepsilon, \gamma$, and the initial data, we proved the global existence and large-time behavior of strong and classical solutions away from vacuum to the  compressible quantum  Navier-Stokes equations with arbitrarily large initial data. This result extends the previous ones on the construction of global strong large-amplitude solutions of the  compressible quantum  Navier-Stokes equations  to the case  $\nu\neq\varepsilon$. Moreover, the vanishing dispersion limit for the classical solutions of the quantum  Navier-Stokes equations is also established  with certain  convergence rates. The proof is based on a new effective velocity which converts the quantum  Navier-Stokes equations into a parabolic system,  and  some elaborate estimates  to derive  the  uniform-in-time positive lower and upper bounds on the specific volume.

\bigbreak
\noindent


\end{abstract}

\section{Introduction }
\setcounter{equation}{0}
The Cauchy problem of the  one-dimensional barotropic compressible quantum  Navier-Stokes equations in the Eulerian coordinates reads as:
\begin{eqnarray}\label{1.1}
\left\{\begin{array}{ll}
    \rho_{\tau}+(\rho u)_{y}=0,\\[2mm]
    (\rho u)_{\tau}+(\rho u^{2}+p(\rho))_{y}=2\nu(\rho u_{y})_{y}+\displaystyle2\varepsilon^{2}\rho\left(\frac{(\sqrt{\rho})_{yy}}{\sqrt{\rho}}\right)_{y},
\end{array}\right.\,\,\,\tau>0, y\in\mathbb{R}
\end{eqnarray}
with the  initial data
\begin{eqnarray}\label{1.2}
&&(\rho(\tau,y),u(\tau,y))|_{\tau=0}=\left(\rho_{0}(y),u_{0}(y)\right)\, \mbox{for}\, y\in\mathbb{R}, \,\,\mbox{and} \inf_{y\in\mathbb{R}}\rho_0(y)>0,\\
&&\lim\limits_{|y|\to+\infty}\left(\rho_{0}(y),u_{0}(y)\right)
=\left(\bar{\rho},\bar{u}\right).
\end{eqnarray}
Here, $\tau$ is the the time variable,  and $y$ is the  space variable.  $\rho=\rho(t,x)>0, u=u(t,x)$ and $p=p(\rho(t,x))$ are unknown functions, which stands for  the density, velocity, and pressure of the fluids,  respectively.  $\nu>0$ is the viscosity constant,   $\varepsilon>0$ is the Planck constant, and $\bar{\rho}>0$, $\bar{u}\in\mathbb{R}$ are given constants.  The expression $\frac{(\sqrt{\rho})_{yy}}{\sqrt{\rho}}$ is called the Bohm potential, which can be interpreted as a quantum potential.

The compressible quantum  Navier-Stokes equations (\ref{1.1}) can be derived from a Chapman-Enskog expansion for the Wigner equation with a BGK term,  see \cite{Brull-derivation,Jungel-derivation} for details.  The system (\ref{1.1}) is also a special case of the general compressible Navier-Stokes-Korteweg equations with density-dependent viscosity and capillarity, which reads as:
\begin{eqnarray}\label{1.3}
\left\{\begin{array}{ll}
    \rho_{\tau}+(\rho u)_{y}=0,\\[2mm]
    (\rho u)_{\tau}+(\rho u^{2}+p(\rho))_{y}=(\mu(\rho) u_{y})_{y}+\displaystyle K_{y},
\end{array}\right.\,\,\,\tau>0, y\in\mathbb{R},
\end{eqnarray}
 where $\mu(\rho)$ is the viscosity coefficient,  and $K$ is the Korteweg tensor given by
 \[K=\rho\kappa(\rho)\rho_{yy}+
          \displaystyle\frac{1}{2}\left(\rho\kappa^\prime(\rho)-\kappa(\rho)\right)\rho_y^2\]
with $\kappa(\rho)$ being  the capillarity  coefficient. Indeed, by choosing $\mu(\rho)=2\nu\rho$ and $\kappa(\rho)=\frac{\varepsilon^2}{\rho}$, then  system (\ref{1.3}) is reduced to the system (\ref{1.1}). The inviscid counterpart of system (\ref{1.1}) is the quantum hydrodynamics model, which has wide applications in quantum fluids,  for instance,  it can be used to describe superfluids \cite{Loffredo-Morato-1993},  quantum semiconductors \cite{Ferry-Zhou-1993}, weakly interacting Bose gases \cite{Grant-1973}, and quantum trajectories of Bohmian mechanics \cite{Wyatt-2005}. The global existence of finite energy weak solutions for the quantum hydrodynamics model has been obtained in \cite{Antonelli-Marcati-ARMA-2012,Antonelli-Marcati-CMP-2009}.

The compressible quantum Navier-Stokes equations has been studied extensively in the last several decades.  Concerning the existence of weak solutions,  J\"{u}ngel \cite{Jungel-SIMA} proved  the  global existence of finite energy weak solutions in the torus $\mathbb{T}^d (d=2,3)$ with the physical parameters $\varepsilon,\nu$, and $\gamma$ satisfying  $\varepsilon>\nu, \gamma\geq1$, if $d=2$, and $\varepsilon>\nu, \gamma>3$ if $d=3$. Then such a result was extended to the cases of $\varepsilon=\nu$ and  $\varepsilon<\nu$ by Dong \cite{Dong} and Jiang \cite{Jiang-NARWA},  respectively. Recently, Antonelli and Spirito \cite{Antonelli-Spirito-ARMA} further improved the result of \cite{Jungel-SIMA} to the cases of $\varepsilon<\nu, \gamma>1$ if $d=2$, and $\varepsilon<\nu<\frac{3\sqrt{2}}{4}\varepsilon, 1<\gamma<3$ if $d=3$, where the authors use a new regular approximating system. We also refer to \cite{Antonelli-Spirito-JHDE,L-Zhang-Zhong-2019,Tang-Zhang-2019} for some improvements of \cite{Antonelli-Spirito-ARMA}. Lacroix-Violet and Vasseur \cite{Lacroix-Violet-Vasseur-JMPA} construct the global finite energy weak solutions which are uniform with respect to the Plank constant $\varepsilon$ for the compressible quantum Navier-Stokes equations in the torus $\mathbb{T}^d$ with $d=2,3$.  The global existence of weak solutions for the compressible quantum Navier-Stokes equations with some additional term was also investigated by some authors, for instance, Gisclon and  Lacroix-Violet \cite{Gisclon-Lacroix-Violet-NA} considered a cold pressure term, and Vasseur and Yu \cite{Vasseur-Yu} added  a damping term.

For the singular limit of weak solutions for the compressible quantum Navier-Stokes equations, Yang et al. \cite{Yang-Ju-Yang} discussed the combined  incompressible limit and semiclassical limit of  global weak solutions  to the local strong solution of the incompressible Navier-Stokes equations with well-prepared initial data in the torus $\mathbb{T}^3$, while the case of ill-prepared initial data was treated in \cite{Kwon-Li}. Just recently, Antonelli, Hientzsch, and Marcati \cite{Antonelli-Hientzsch-Marcati-SIMA} considered the incompressible limit of global weak solutions for the quantum Navier-Stokes equations in the whole space $\mathbb{R}^3$, where, in contrast to \cite{Yang-Ju-Yang,Kwon-Li}, the global weak solution of the incompressible Navier-Stokes equations is obtained as the limit solution. In \cite{Lacroix-Violet-Vasseur-JMPA, Gisclon-Lacroix-Violet-NA}, the authors studied the semiclassical limit of global weak solutions of the quantum Navier-Stokes equations  to the global weak solutions of the compressible Navier-Stokes equations in the torus $\mathbb{T}^d$ with $d=2,3$. Note that the initial data in the above references \cite{Jungel-SIMA,Jiang-NARWA,Dong,Antonelli-Spirito-JHDE,L-Zhang-Zhong-2019,Tang-Zhang-2019,Gisclon-Lacroix-Violet-NA,Vasseur-Yu,Antonelli-Spirito-ARMA,Antonelli-Hientzsch-Marcati-SIMA,Lacroix-Violet-Vasseur-JMPA,Yang-Ju-Yang,Kwon-Li} concerning the weak solution of the compressible quantum Navier-Stokes equations can be arbitrarily large.

For the global strong solutions to the compressible quantum Navier-Stokes equations with large initial data,  relatively fewer results are obtained. DiPerna \cite{DiPerna}, J\"{u}ngel \cite{Jungel-NA}, and Charve and Haspot \cite{F-Charve-B-Haspot-2011}  proved  the global existence of strong  solutions with large initial data for the one-dimensional compressible quantum Navier-Stokes equations (\ref{1.1}) in the whole space $\mathbb{R}$ for the case $\varepsilon=\nu$. Recently, the large time behavior of strong solutions to the Cauchy problem  (\ref{1.1})-(\ref{1.2})  with large initial data and  $\varepsilon=\nu$  was obtained by Chen and Li \cite{Chen-Li-2021}.  Haspot \cite{Haspot-2017-MA} studied the global existence of strong solution for the compressible quantum Navier-Stokes equations  with large  initial data and $\varepsilon=\nu$ in the multi-dimensional whole space $\mathbb{R}^N$ with $N\geq2$.

However, as far as we know,  there is no result on the global existence of strong solutions to the compressible quantum Navier-Stokes equations with large initial data for the case $\varepsilon\neq\nu$ up to now. The aim of this paper is to study this problem, and as a first step, we shall consider the global existence and vanishing dispersion limit of  strong/classical solutions for the Cauchy problem  of the compressible quantum Navier-Stokes equations (\ref{1.1}).

Throughout this paper,  we consider the $\gamma-$law pressure, i.e,  $p(\rho)=\rho^{\gamma}$ with $\gamma\geq1$ being a constant, and the far-fields of the initial data $(\bar{\rho}, \bar{u})$ is assumed to be $(1,0)$ for simplicity.
Suppose that $(\rho(\tau,y),u(\tau,y))$ is a classical solution of the Cauchy problem (\ref{1.1})-(\ref{1.2}), then one has
\begin{equation}\label{1.4}
2\varepsilon^{2}\rho\left(\frac{(\sqrt{\rho})_{yy}}{\sqrt{\rho}}\right)_{y}=\varepsilon^{2}\left(\rho_{yy}-\frac{\rho_{y}^{2}}{\rho}\right)_{y}.
\end{equation}
 Setting  the coordinate transformation:
$$x=\int_{0}^{y}\rho(\tau,z)dz-\int_{0}^{\tau}\rho(\xi,0)u(\xi,0)d\xi, \quad   t=\tau,$$
and noticing (\ref{1.4}), then the Cauchy problem (\ref{1.1})-(\ref{1.2}) in the Lagrangian coordinates becomes
\begin{eqnarray}\label{1.5}
\left\{\begin{array}{ll}
    v_{t}-u_{x}=0,\\[2mm]
    u_{t}+p(v)_{x}=\displaystyle2\nu\left(\frac{u_{x}}{v^{2}}\right)_{x}
    +\varepsilon^{2}\left(-\frac{v_{xx}}{v^{4}}+\frac{2v_{x}^{2}}{v^{5}}\right)_{x},
\end{array}\right.\,\,\,x\in\mathbb{R},t>0
\end{eqnarray}
with the initial data
\begin{eqnarray}\label{1.6}
(v(t,x),u(t,x))|_{t=0}=\left(v_{0}(x),u_{0}(x)\right),\quad \lim\limits_{|x|\to+\infty}\left(v_{0}(x),u_{0}(x)\right)=(1,0),
\end{eqnarray}
where $v=\frac{1}{\rho}>0$ denotes the specific volume. The term $\displaystyle\varepsilon^2\left(-\frac{v_{xx}}{v^{4}}+\frac{2v_{x}^{2}}{v^{5}}\right)_{x}$ represents the dispersive effect of the system, and the  regimes $\varepsilon\ll\nu$, $\varepsilon\sim\nu$ and $\varepsilon\gg\nu$ are called {\it the  parabolic regime, the intermediary regime} and {\it the dispersive regime}, respectively. In the parabolic regime, the viscosity dominates so the parabolic effect is predominant, and in  the dispersive regime,   the capillarity dominates so the dispersive effect is predominant.

To state our main results, let us first introduce the notations used throughout this paper. The symbol $C$ stands for some generic positive constant which may vary from line to line, and  $C(\cdot,\cdots,\cdot)$ or  $C_i(\cdot,\cdots,\cdot)(i\in {\mathbb{N}})$ denotes the constant which depends explicitly on the  qualities listed in brackets. $C^k_b(\mathbb{R})(k\geq0)$ is  the space of the $k$-times continuously differentiable functions in $\mathbb{R}$ whose derivatives up to  $k$-th order are all bounded. $L^p(\mathbb{R})(1\leq p\leq+\infty)$ is  the usual Lebesgue space with the norms  $\|\cdot\|_{L^p}:=\|\cdot\|_{L^p(\mathbb{R})}=\left(\int_{\mathbb{R}}|\cdot|^pdx\right)^{\frac{1}{p}}$ for $p=2$, and  $\|\cdot\|_{L^\infty}:=\|\cdot\|_{L^\infty(\mathbb{R})}=ess sup_x|\cdot|$ for $p=\infty$.   For a non-negative real number $s$, $H^s(\mathbb{R})$ is the standard  Sobolev space with its norm $\|f\|_{s}:=\|f\|_{H^{s}(\mathbb{R})}=(\int_{\mathbb{R}}(1+\xi^2)^{s}|\hat{f}(\xi)|^2d\xi)^{\frac{1}{2}}$, where $\hat{f}(\xi)$ is the fourier transformation of the function $f(x)$. When $s$ is a non-negative integer,  the norm $\|f\|_{s}$ is equivalent to  $\left(\sum_{i=0}^s\|\partial_x^i f\|_{L^2}^2\right)^{\frac{1}{p}}$. For brevity,  we denote $\|\cdot\|:=\|\cdot\|_{L^2}$,  and  $\|\cdot\|_{L_{T,x}^\infty}:=\|\cdot\|_{L^\infty([0,T]\times\mathbb{R})}$.

 Our first  result is  on the global existence of strong solution for the Cauchy problem (\ref{1.5})-(\ref{1.6}) when  the parameters $(\varepsilon,\nu)$ lie in  the intermediary regime (i.e., $\varepsilon\sim\nu$),  which can be stated as follows.

\begin{Theorem}
Let $0<\varepsilon\leq\nu$ and $\gamma\geq\frac{7}{4}$. Suppose that the initial data
 $(v_0(x),u_0(x))$  satisfy $(v_0(x)-1,u_0(x))\in H^2(\mathbb{R})\times H^1(\mathbb{R})$,  and  $\underline{V}\leq v_0(x)\leq\overline{V}$ for all $x\in\mathbb{R}$ and some given positive constants $\underline{V},\overline{V}$.
Then the Cauchy problem (\ref{1.5})-(\ref{1.6}) admits a unique global-in-time  strong solution $(v(t,x), u(t,x))$  satisfying
\begin{equation}\label{1.7}
C_0^{-1}\leq v(t,x)\leq C_0,\,\, \forall\,(t,x)\in[0,\infty)\times\mathbb{R},
\end{equation}
\begin{eqnarray}\label{1.8}
&&\|v(t)-1\|^2_{2}+\|u(t)\|^2_{1}+\int_0^t\left(\|v_x(\tau)\|^2_{2}+\|u_x(\tau)\|^2_{1}\right)d\tau\nonumber\\
&&\leq C_1\left(\|v_0-1\|^2_{2}+\|u_0\|^2_{1}\right), \,\, \forall\, t>0.
\end{eqnarray}
Moreover, $(v(t,x), u(t,x))$ satisfies the following large-time behavior:
\begin{equation}\label{1.9}
\lim_{t\rightarrow+\infty}\left\{\left\|\left(v(t)-1,u(t)\right)\right\|_{L^\infty}+\left\|v_x(t)\right\|_{1}+\left\|u_x(t)\right\|\right\}=0.
\end{equation}
Here $C_0$ and $C_1$ are two positive constants depending only on $\varepsilon, \nu, \gamma, \underline{V}, \overline{V}, \|v_0-1\|_{2}$, and $\|u_0\|_{1}$.
\end{Theorem}

If the initial data $(v_0(x),u_0(x))$ is sufficiently regular and the  parameters $(\varepsilon,\nu)$ lie in  the parabolic regime (i.e., $\varepsilon$ is sufficiently small), then we can obtain the global existence of classical solutions to the Cauchy problem (\ref{1.5})-(\ref{1.6}). The precise statement of this result is as follows.

\begin{Theorem}
Let  $\displaystyle1\leq\gamma\leq2$ and $\nu$ be any fixed positive constant. Suppose that $(v_0(x)-1,u_0(x))\in H^5(\mathbb{R})\times H^4(\mathbb{R})$,  and $\underline{V}\leq v_0(x)\leq\overline{V}$ for all $x\in\mathbb{R}$ and some given positive constants $\underline{V}, \overline{V}$.
Then there exists a sufficiently small positive constant $\varepsilon_0$ depending only on $\nu,\gamma,  \underline{V}, \overline{V}$, $\|v_0-1\|_{5}$, and $\|u_0\|_4$,  such that if $0<\varepsilon\leq\varepsilon_0$, the Cauchy problem (\ref{1.5})-(\ref{1.6}) has a unique global-in-time  classical solution $(v(t,x), u(t,x))$ satisfying 
\begin{equation}\label{1.13}
C_2^{-1}\leq v(t,x)\leq C_2,\,\, \forall\,(t,x)\in[0,\infty)\times\mathbb{R},
\end{equation}
\begin{eqnarray}\label{1.14}
&&\|v(t)-1\|^2_{5}+\|u(t)\|^2_{4}+\varepsilon^2\|v_x(t)\|^2_{4}+
\int_0^t\left(\|(u_x(\tau),v_x(\tau))\|^2_{4}+\varepsilon^2\|v_{xx}(\tau)\|^2_{4}\right)d\tau\nonumber\\
&&\leq C_3\left(\|v_0-1\|^2_{5}+\|u_0\|^2_{4}\right), \,\,\forall\, t>0,
\end{eqnarray}
and
\begin{equation}\label{1.15}
\lim_{t\rightarrow+\infty}\left\{\left\|\left(v(t)-1,u(t)\right)\right\|_{L^\infty}+\left\|v_x(t)\right\|_{4}+\left\|u_x(t)\right\|_{3}\right\}=0.
\end{equation}
Here $C_2$ is a  positive constant depending only on $\nu, \gamma,\underline{V}, \overline{V},  \|v_0-1\|_{1}$, and $\|u_0\|$, and $C_3$ is a  positive constant depending only on $\nu, \gamma,\underline{V}, \overline{V},  \|v_0-1\|_{5}$, and $\|u_0\|_4$.
\end{Theorem}
\begin{Remark}We list some remarks  on  Theorems 1.1-1.2 as follows.
\begin{itemize}
\item[$(1)$] In Theorems 1.1-1.2, the uniform-in-time positive lower and upper bounds of the specific volume $v(t,x)$ are deduced by using different methods, and in  these procedures, we need to ask  that the initial data $(v_0(x)-1,u_0(x))$ belong to $H^2(\mathbb{R})\times H^1(\mathbb{R})$ and $H^5(\mathbb{R})\times H^4(\mathbb{R})$, respectively. Note that there are no smallness assumptions on the initial data $(v_0(x)-1,u_0(x))$.

\item[$(2)$] In Theorems 1.1-1.2, the parameters $\varepsilon$ and $\nu$ need to satisfy $0<\varepsilon\leq\nu$ or $0<\varepsilon\ll\nu$, i.e., only the parabolic regime and some subdomains of the intermediary regime are covered here. The case of $\varepsilon>\nu>0$ is not discussed in the present paper, which includes the dispersive regime (i.e., $\varepsilon\gg\nu>0$). When $\varepsilon>\nu$, if we follow the same approach as that in the proof of Theorem 1.1,  then the constants $c_+$ and $c_-$ in the system (\ref{2.22}) should be replaced by \[c_+=\nu+i\sqrt{\varepsilon^2-\nu^{2}}, \quad c_-=\nu-i\sqrt{\varepsilon^2-\nu^{2}}.\]
    Consequently, the system (\ref{2.22}) becomes a complex one, for which our analysis for  the derivation of the $L^\infty$ bound on the effective velocity $\xi(t,x)=u(t,x)-c_+\frac{v_x(t,x)}{v^2(t,x)}$ can not be applied (see the proof of Lemma 2.6). Indeed, in our analysis, we need to assume that the constants $c_+$ and $c_-$ in (\ref{2.25}) are  positive. That is why we only consider the case of $\varepsilon\leq\nu$ throughout this paper. Theorems 1.1-1.2 extend the previous results  \cite{Jungel-NA, F-Charve-B-Haspot-2011, DiPerna, Chen-Li-2021, Haspot-2017-MA} on the construction of global strong large-amplitude solutions of the compressible quantum Navier-Stokes equations to the case $\nu\neq\varepsilon$.

\item[$(3)$] The result of Theorem 1.1 also holds for the case of $\nu=\varepsilon$, which has been studied in \cite{Jungel-NA, F-Charve-B-Haspot-2011, Chen-Li-2021}. Indeed,  the authors in \cite{Jungel-NA, F-Charve-B-Haspot-2011, Chen-Li-2021} transformed the system  (\ref{1.1}) with $\nu=\varepsilon$ into a parabolic one via the effective velocity $w(t,x)=u(t,x)+\frac{\nu\rho_x(t,x)}{\rho(t,x)}$, then the transformed system has positively invariant region for the unknown functions $(\rho(t,x), w(t,x))$, which yields the $L^\infty$ bounds for $\rho(t,x)$ and $w(t,x)$ immediately (such a result was first obtained by DiPerna in \cite{DiPerna} for the one-dimensional compressible Navier-Stokes equations with artificial viscosity). In this paper, we derive the $L^\infty$ bounds on $\rho(t,x)$ and $w(t,x)$ by the pure energy method rather than the theory of positively invariant region.  And such an energy method can also be applied to study the case of $\nu\neq\varepsilon$ (see the proof of Theorem 1.1 for details). However, we need the assumption that $\gamma\geq\frac{7}{4}$ rather than $\gamma\geq1$ in \cite{DiPerna, Jungel-NA, F-Charve-B-Haspot-2011, Chen-Li-2021}.

\item[$(4)$] The global existence of strong large solutions to the multidimensional barotropic  compressible quantum Navier-Stokes equations with $\varepsilon\neq\nu$, as well as the nonisothermal compressible quantum Navier-Stokes equations are  still  unknown up to now. These interesting problems will be pursued by the authors in the future.
\end{itemize}
\end{Remark}

Finally, we consider the vanishing dispersion limit for the one-dimensional compressible quantum  Navier-Stokes equations (\ref{1.5}), which  corresponds to the physical problem of letting the Planck constant $\varepsilon$ tend to zero, while keeping the viscosity coefficient $\nu$ be a constant. Thanks to the uniform-in-$\varepsilon$ estimate (\ref{1.14}), we can obtain the global classical solutions for the one-dimensional compressible Navier-Stokes equations with large data through the  vanishing dispersion limit (i.e., letting $\varepsilon\rightarrow0$). In order to  state this result precisely, we denote  the classical solutions to the Cauchy problem (\ref{1.5})-(\ref{1.6}) with the  superscripts $\varepsilon$ as $\left(v^{\varepsilon}(t,x), u^{\varepsilon}(t,x)\right)$, then formally as $\varepsilon\rightarrow0$, we see that the limit solution $\left(v^{0}(t,x), u^{0}(t,x)\right)$ satisfies the following compressible Navier-Stokes equations:
\begin{eqnarray}\label{1.16}
\left\{\begin{array}{ll}
    v^0_{t}-u^0_{x}=0,\\[2mm]
    u^0_{t}+p(v^0)_{x}=\displaystyle2\nu\left(\frac{u^0_{x}}{(v^{0})^2}\right)_{x},
\end{array}\right.\,\,\,x\in\mathbb{R},t>0
\end{eqnarray}
with initial data
\begin{eqnarray}\label{1.17}
\left(v^{0}(t,x),u^{0}(t,x)\right)|_{t=0}=\left(v_{0}(x),u_{0}(x)\right),\,\,\,x\in\mathbb{R}.
\end{eqnarray}
Then our last main result on the vanishing dispersion limit $\varepsilon\rightarrow0$ is stated as follows.
\begin{Theorem}
Let  $(v^{\varepsilon}(t,x), u^{\varepsilon}(t,x))$ be the global classical solution of the Cauchy problem (\ref{1.5})-(\ref{1.6}) obtained in Theorem 1.2. Then as $\varepsilon\rightarrow0$, it holds that
\[v^\varepsilon(t,x)\rightarrow v^0(t,x) \quad \mbox{strongly in} \,\,C\left(0,T; H^{5-s}_{loc}(\mathbb{R})\bigcap C_b^4(\mathbb{R})\right),\]
\[u^\varepsilon(t,x)\rightarrow u^0(t,x) \quad \mbox{strongly in} \,\,C\left(0,T; H^{4-s}_{loc}(\mathbb{R})\bigcap C_b^3(\mathbb{R})\right)\]
for any $s\in(0,\frac{1}{2})$ and any fixed $T>0$, where $(v^{0}(t,x),u^{0}(t,x))$ is the global classical solution to the Cauchy problem (\ref{1.16})-(\ref{1.17}),  which satisfies
\begin{eqnarray}\label{1.18-1}
C_2^{-1}\leq v^0(t,x)\leq C_2,\,\, \forall\,(t,x)\in[0,\infty)\times\mathbb{R},
\end{eqnarray}
\begin{eqnarray}\label{1.18}
\left\|v^0(t)-1\right\|^2_{5}+\|u^0(t)\|^2_{4}+
\int_0^t\left\|\left(u^0_x(\tau),v^0_x(\tau)\right)\right\|^2_{4}\,d\tau\leq C_3\left(\|v_0-1\|^2_{5}+\|u_0\|^2_{4}\right), \,\, \forall\, t>0,
\end{eqnarray}
and
\begin{equation}\label{1.18-2}
\lim_{t\rightarrow+\infty}\left\{\left\|\left(v^0(t)-1,u^0(t)\right)\right\|_{L^\infty} +\left\|\left(v^0_x(t),u_x^0(t)\right)\right\|_{3}\right\}=0,
\end{equation}
where $C_2$ and $C_3$ are positive constants given in Theorem 1.2.

Moreover, it holds for any fixed time $t\in[0,+\infty)$ that
\begin{eqnarray}\label{1.19}
\left\|\left(v^\varepsilon(t)-v^0(t),u^\varepsilon(t)-u^0(t)\right)\right\|&\leq& C_4\varepsilon^2t^{\frac{1}{2}}e^{\lambda t},\\[2mm]
\left\|\left(\frac{\partial v^\varepsilon(t)}{\partial x^k}-\frac{\partial v^0(t)}{\partial x^k},\frac{\partial u^\varepsilon(t)}{\partial x^k}-\frac{\partial u^0(t)}{\partial x^k}\right)\right\|&\leq& C_5\varepsilon^2 \left(1+t^{\frac{k+1}{2}}\right)e^{\lambda t},\,\, k=1,2,3,\label{1.20}
\end{eqnarray}
where $C_4,C_5$, and $\lambda$ are positive constants depending only on $\nu, \gamma,\underline{V}, \overline{V},  \|v_0-1\|_{5}$, and $\|u_0\|_4$.
\end{Theorem}

\begin{Remark}Two  remarks  on  Theorem 1.3 are  presented  as follows.
\begin{itemize}
\item[$(1)$]Compared to (\ref{1.8}),  the estimate (\ref{1.14}) is uniform with respect to the Plank constant $\varepsilon$, and from which, we can discuss the vanishing dispersion limit $\varepsilon\rightarrow0$ in Theorem 1.3, which establishes  the relationship between the compressible  quantum  Navier-Stokes equations (\ref{1.5}) and the compressible Navier-Stokes equations (\ref{1.16}) in the case of global classical solutions with large data. Furthermore, our result on the global large solutions to the compressible Navier-Stokes equations (\ref{1.16}) can also cover the case of Saint-Venant system for the shallow water, i.e., the case of $\gamma=2$.

\item[$(2)$] The vanishing dispersion limit $\varepsilon\rightarrow0$ for the  compressible  quantum  Navier-Stokes equations in the torus $\mathbb{T}^d (d=2,3)$ was  discussed by Lacroix-Violet and Vasseur \cite{Lacroix-Violet-Vasseur-JMPA}, and  Gisclon and Lacroix-Violet \cite{Gisclon-Lacroix-Violet-NA},  where they obtained the global weak solutions of the barotropic compressible Navier-Stokes equations with large data as the limit solutions.  We also mention that the vanishing  capillarity limit for the Cauchy problem of the barotropic   compressible Navier-Stokes-Korteweg equations was studied by Bian, Yao and Zhu \cite{BYZ}, and Burtea and Haspot\cite{Burtea-Haspot-2022},  where the global smooth small  solution of the three-dimensional compressible Navier-Stokes equations and the global strong large solution of the one-dimensional compressible Navier-Stokes equations are obtained respectively, as the capillarity coefficient tends to zero. In the above references,  the large-time behavior of the limit solutions are not shown, and the convergence rate estimates are only given  in \cite{BYZ}. Moreover,  the authors in \cite{BYZ} studied the Navier-Stokes-Korteweg system with constant viscosity and capillarity, while the case of density-dependent viscosity and capillarity is treated in \cite{Burtea-Haspot-2022}.
   A typical example of the  viscosity coefficient $\mu(\rho)$ and the capillarity  coefficient $\kappa(\rho)$ imposed in \cite{Burtea-Haspot-2022} is $\mu(\rho)=\rho^\alpha$  and $\kappa(\rho)=\rho^{2\alpha-3}$ with $\alpha\in(0,\frac{1}{2})$, which  can not cover the case of quantum Navier-Stokes equations considered in the present paper (i.e., the case of $\alpha=1$).
   Compared to \cite{Lacroix-Violet-Vasseur-JMPA,Gisclon-Lacroix-Violet-NA,BYZ,Burtea-Haspot-2022}, we obtain not only the  global existence but also the large-time behavior and convergence rate of the  smooth large solution of the one-dimensional barotropic compressible  Navier-Stokes equations in Theorem 1.3 by making full use of the one-dimensional feature.
\end{itemize}
\end{Remark}

Now we outline the main ideas used in proving Theorems 1.1-1.3. As pointed out in \cite{Jungel-NA, F-Charve-B-Haspot-2011, DiPerna, Chen-Li-2021, Haspot-2017-MA}, the key ingredient in establishing the global existence of nonvacuum strong/classical solutions to the one-dimensional compressible quantum Navier-Stokes equations (\ref{1.1}) is to deduce  the uniform-in-time positive lower and upper bounds on the specific volume $v(t,x)$. However, the methods used in \cite{Jungel-NA, F-Charve-B-Haspot-2011, DiPerna, Chen-Li-2021, Haspot-2017-MA} for the barotropic compressible quantum Navier-Stokes equations with $\nu=\varepsilon$ can not be applied to our problem (\ref{1.5})-(\ref{1.6}). Indeed, as mentioned in Remark 1.1 (3), if one follows the methods of \cite{Jungel-NA, F-Charve-B-Haspot-2011, DiPerna, Chen-Li-2021}, that is, transform the system (\ref{1.1}) through the effective velocity $w(t,x)=u(t,x)+\frac{\nu\rho_x(t,x)}{\rho(t,x)}$, then the transformed system does not have any positively invariant region due to the general assumption that $\nu\neq\varepsilon$. Consequently, the theory of positively invariant region used in \cite{DiPerna, Jungel-NA, F-Charve-B-Haspot-2011, Chen-Li-2021} can not be applied here any more.
Moreover, the analysis in \cite{Haspot-2017-MA} for the mutil-dimensional compressible quantum Navier-Stokes equations relies crucially on the assumption that the pressure $p(\rho)=a\rho$ with the constant $a>0$, which can not be applied to the general $\gamma$ pressure law with $\gamma>1$. To achieve this key point, some new mathematical techniques and useful estimates are developed  in  the proof of Theorems 1.1-1.2, which can be stated  as follows:
\begin{itemize}
\item[$\bullet$] Firstly, in the proof of Theorem 1.1, the uniform-in-time lower bound of $v(t,x)$ can be obtained easily from the basic energy estimates (\ref{2.4}) and the H\"{o}lder inequality by using the idea developed by Jiang, Xin, and Zhang in \cite{Jiang-Xin-Zhang} (see also \cite{Guo-Jiu-Xin,Li-Li-Xin}), the details can be found in the proof of Lemma 2.3 below. Then in order to deduce the  upper  bound of  $v(t,x)$, we perform the estimates on $\int_{\mathbb{R}}\frac{\nu v_x^2}{v^4}dx$,  $\varepsilon^2\int_0^t\int_{\mathbb{R}}\frac{v_{xx}^2}{v^6}dxd\tau$, and $\varepsilon^2\int_0^t\int_{\mathbb{R}}\frac{v_{x}^4}{v^8}dxd\tau$, which rely  on the interesting Germain and LeFloch's inequality introduced in \cite{Germain-LeFloch-2012} (see Lemma 2.4). Based on these, we obtain a novel  $L^2(0,T; L^\infty(\mathbb{R}))$ estimate  on the function $u(t,x)/v^\frac{3}{4}(t,x)$, which plays an important role in establishing  the $L^\infty$ bound of the new effective velocity  $\xi(t,x):=u(t,x)-c_+\frac{v_x(t,x)}{v^2(t,x)}$ with $c_+=\nu+\sqrt{\nu^2-\varepsilon^2}$. Notice  that this  new effective velocity  $\xi(t,x)$ will be reduced to  $\xi(t,x)=u(t,x)-\nu\frac{v_x(t,x)}{v^2(t,x)}$ if $\nu=\varepsilon$, which is exactly  the effective velocity $w(t,x)=u(t,x)+\frac{\nu\rho_x(t,x)}{\rho(t,x)}$ in the Eulerian  coordinates. The assumption $\gamma\geq\frac{7}{4}$ is also employed  in  the proof of $\|\xi(t)\|_{L^\infty}$ (see the proof of Lemmas 2.5 and 2.6 for details). With the aid of the effective velocity  $\xi(t,x)$,  the original system (\ref{1.5}) turns into a parabolic one for the unknown functions $(v(t,x), \xi(t,x))$. Then from the parabolic equation of the specific volume $(\ref{2.22})_1$, we can obtain the estimate of $\|\phi(t)\|_{L^4}$ with $\phi(t,x)=v(t,x)-1$, which together with the H\"{o}lder inequality yields the upper bound of $v(t,x)$ (see Lemma 2.7 and its proof, here the upper bound of $v(t,x)$  depends on the time step $T$). In order to show that the upper bound of $v(t,x)$ is uniform in time, we use a large-time stability analysis which is similar to  that of Lemma 4.3 in  \cite{Chen-Li-2021}.
\item[$\bullet$] Secondly, since the lower and upper bounds of $v(t,x)$ in Theorem 1.2 should be independent of  the Planck constant $\varepsilon$, we derive the uniform-in-time lower bound of $v(t,x)$ from the estimate of $\int_{\mathbb{R}}\frac{\nu v_x^2}{v^4}dx$ instead of the basic energy estimates (\ref{2.4}). Then to show the upper bound of $v(t,x)$, we introduce another effective velocity $\omega(t,x):=u(t,x)-2\nu\frac{v_x(t,x)}{v^2(t,x)}$, which makes the specific volume $v(t,x)$ satisfy a parabolic equation, while the effective velocity $\omega(t,x)$ satisfy a first order nonlinear equation with a damping  term (see (\ref{3.8})). Then we can obtain a  pointwise  $L^\infty$ estimate  of the specific volume $v(t,x)$  by modifying the argument from the proof of Theorem 2.1 in \cite{Burtea-Haspot}, where the upper bound of the specific volume $v(t,x)$ is obtained for the Cauchy problem of the one-dimensional compressible Navier-Stokes equations. Furthermore, we can deduce the energy estimates of the solutions based on the a priori assumption that $v(t,x)\leq M$ for all $(t,x)\in[0,T]\times\mathbb{R}$ and some positive constant $M>0$ (see (\ref{3.3})-(\ref{3.4}) for details).  Then by a delicate continuation argument, we can obtain the uniform-in-time lower and upper bounds of $v(t,x)$.
\end{itemize}

Once the uniform-in-time positive lower and upper bounds of the specific volume $v(t,x)$ is obtained, the uniform-in-time energy estimates of solutions to the Cauchy problem (\ref{1.5})-(\ref{1.6}) can be derived  easily, and finally,  Theorems 1.1-1.2 can be proved by the standard continuation argument. Here we would like to point out that  in deducing the upper bounds of $v(t,x)$ for Theorem 1.2,  one needs to deal with the dispersion term $\varepsilon^{2}\left(-\frac{v_{xx}}{v^{4}}+\frac{2v_{x}^{2}}{v^{5}}\right)_{x}$ in the energy type estimates. To overcome this difficulty,  we make the a priori assumption (\ref{3.1}), and use the smallness assumption on the Plank constant $\varepsilon$ (see (\ref{3.20}) for details).

The proof of Theorem 1.3 is motivated by \cite{BYZ} for the vanishing capillarity limit of smooth small-amplitude solutions to the three-dimensional isentropic compressible Navier-Stokes-Korteweg equations. In the present paper, the limit as the Plank constant $\varepsilon\rightarrow0$ is proved by employing the Lions-Aubin lemma. And by the energy method, the convergence rates of solutions of the compressible quantum Navier-Stokes equations (\ref{1.5}) to the solutions of the  compressible  Navier-Stokes equations (\ref{1.16}) are also obtained.

Before finishing this section,  we point out that there have been many mathematical results on the compressible Navier-Stokes equations with density-dependent viscosity. We refer to \cite{Burtea-Haspot,Constantin-2019,Chen-Zhao-2020,Haspot-2018,Mellet-Vasseur-2008,Qin-Yao} and the references therein for the global existence and large-time behavior of strong solutions for the one-dimensional compressible Navier-Stokes equations without vacuum,\cite{Li-Li-Xin,Liu-Xin-Yang-1998,Jiang-Xin-Zhang,Jiu-Xin-2008,Yang-Yao-Zhu,Zhang-Fang,Yang-Zhao-2002,Zhu-2010,Jiu-Wang-Xin-2011} and the references therein for the local/global existence and large-time behavior of weak solutions for the one-dimensional compressible isentropic Navier-Stokes equations with vacuum,\cite{Guo-Jiu-Xin,Jiu-wang-xin,Huang-Li,Li-Pan-Zhu} and the references therein for the global existence of weak or strong solutions for the compressible isentropic Navier-Stokes equations in the multi-dimensional case. In this paper, we borrow some ideas from \cite{Jiang-Xin-Zhang,Li-Li-Xin,Guo-Jiu-Xin,Burtea-Haspot}  to show the lower and upper bounds of the specific volume for the one-dimensional quantum Navier-Stokes equations (\ref{1.5}), see the proof of Theorems 1.1-1.2 in the following Sections 2 and 3 for details.

The rest
 of this article unfolds as follows:  Section 2 is devoted to proving Theorem 1.1, where the   novel  $L^2(0,T; L^\infty(\mathbb{R}))$ estimate on the function $u(t,x)/v^\frac{3}{4}(t,x)$ and the $L^\infty$ estimate on  the new effective velocity $\xi(t,x)$ are derived.
In Section 3, we present the proof of Theorem 1.2. Some uniform estimates with respect to the Plank constant $\varepsilon$ for the  strong solutions of  the Cauchy problem (\ref{1.5})-(\ref{1.6}) are given in this section. In Section 4, we study the vanishing dispersion limit as $\varepsilon\rightarrow0$. The convergence rates of solutions to the problem (\ref{1.5})-(\ref{1.6}) is also obtained in this section.  In the final Section 5,  we first list some useful lemmas which will play important roles in our analysis,  and then  give the rigorous derivation of the system (\ref{2.22}).

\section{Global existence in the intermediary regime}
\setcounter{equation}{0}
This section is devoted to proving Theorem 1.1, which is a global existence result in the intermediary regime $\varepsilon\sim\nu$ with large initial data.  First,   we define  the following sets of functions:
\begin{eqnarray}\label{2.1}
&&X_i(0,T;m,M)\nonumber\\
&=&\left\{(v(t,x), u(t,x))\left|
\begin{array}{c}
v(t,x)-1\in C(0, T; H^{i+1}(\mathbb{R})),\,\,u(t,x)\in C(0, T; H^i(\mathbb{R})),\\[2mm]
(v_x(t,x),u_x(t,x))\in L^2(0, T; H^{i+1}(\mathbb{R})\times H^{i}(\mathbb{R})),\\[2mm]  \displaystyle m\leq v(t,x)\leq M, \,\,(t,x)\in[0,T]\times\mathbb{R}
\end{array}
\right.\right\}
\end{eqnarray}
with $i\in\{1,4\}$, and $m, M$, and $T$ being some positive constants.

Then the local existence of solutions to the Cauchy  problem (\ref{1.5})-(\ref{1.6}) is stated as follows. 
\begin{Proposition} [Local existence]
Assume that the conditions  of Theorem 1.1 hold, and the initial data $(v_0(x)-1,u_0(x))\in H^{2}(\mathbb{R})\times H^1(\mathbb{R})$, then there exists a sufficiently small positive constant $t_0=t_0(\underline{V},N_0)$ depending only on $\nu,\varepsilon, \gamma,\underline{V}$, and $N_0:=(\|v_0-1\|^2_{2}+\|u_0\|^2_1)^{\frac{1}{2}}$, such that the Cauchy  problem (\ref{1.5})-(\ref{1.6}) admits a unique  solution $(v(t,x), u(t,x))\in X_1(0,t_0; \frac{1}{2}\underline{V}, 2\overline{V})$, and
 \[
\sup_{[0,t_0]}\{\|v(t)-1\|^2_{2}+\|u(t)\|^2_{1}\}+\int_0^{t_0}\left(\|v_x(\tau)\|^2_{2}+\|u_x(\tau)\|^2_{1}\right)d\tau\nonumber\\
\leq bN_0^2,
\]
 where $b>1$ is a positive constant depending only on $\underline{V}$, $\overline{V}$.
\end{Proposition}

Proposition 2.1 can be obtained by using the dual argument and iteration technique, which is similar to the proof of Theorem  1.1 in \cite{Hattori-Li-1994}.  We omit it here for brevity.

 To show  the global existence of solutions to  the problem (\ref{1.5})-(\ref{1.6}), one needs to derive  certain {\it a priori} estimates on the solution $(v(t,x),u(t,x))$.
\begin{Proposition} [{\it A priori} estimates]
Let  the conditions of  Theorem 1.1 hold, and $(v(t,x), u(t,x))\in X_1(0,T;m,M)$ be a solution of the Cauchy problem  (\ref{1.5})-(\ref{1.6}) defined in  $\Pi_T=[0,T]\times\mathbb{R}$ for some positive constants $m, M$, and $T$.  Then  the following estimates  hold:
\begin{equation}\label{2.2}
C_0^{-1}\leq v(t,x)\leq C_0,\quad \forall\,(t,x)\in[0,T]\times\mathbb{R},
\end{equation}
\begin{eqnarray}\label{2.3}
&&\|v(t)-1\|^2_{2}+\|u(t)\|^2_{1}+\int_0^t\left(\|v_x(\tau)\|^2_{2}+\|u_x(\tau)\|^2_{1}\right)d\tau\nonumber\\
&&\leq C_1\left(\|v_0-1\|^2_{2}+\|u_0\|^2_{1}\right), \quad \forall\, t\in[0,T],
\end{eqnarray}
where $C_0$ and $C_1$ are two  positive constants depending only on $\varepsilon, \nu, \gamma,\underline{V}, \overline{V}, \|v_0-1\|_{2}$, and $\|u_0\|_{1}$.
\end{Proposition}

As usual, Theorem 1.1 will be obtained by the continuation argument which is based on the local existence result Proposition 2.1 and the {\it a priori} estimates stated in Proposition 2.2. Thus it remains to show  Proposition 2.2 holds, which is left for the following subsection.

\subsection{{\it A priori} estimates}
The main target of this subsection is to deduce  the {\it a priori} estimates (\ref{2.2})-(\ref{2.3}).
First of all, we derive the following elementary  energy estimates.
\begin{Lemma}Assume that the conditions  of Proposition 2.2 hold, then  for all $t\in[0,T]$,
\begin{eqnarray}\label{2.4}
 \int_{\mathbb{R}}\left(\Phi(v)+\frac{u^2}{2}\right)dx +\frac{\varepsilon^{2}}{2}\int_{\mathbb{R}}\frac{v_x^2}{v^{4}}dx+2\nu\int_0^t\int_{\mathbb{R}}\frac{ u_x^2}{v^{2}}dxd\tau
 =\int_{\mathbb{R}}\left(\Phi(v_{0})+\frac{u_{0}^2}{2}+\frac{\varepsilon^{2}v_{0x}^2}{2v_{0}^{4}}\right)dx,
\end{eqnarray}
where the function $\Phi(v)$ is defined by
\[\Phi(v):=p(1)(v-1)-\displaystyle\int_{1}^v p(s)ds=\left\{\begin{array}{ll}
 v-1+\displaystyle\frac{1}{\gamma-1}\left(v^{-\gamma+1}-1\right),\quad if\,\, \gamma>1,\\[2mm]
    v-1-\ln v,\quad if\,\, \gamma=1.
\end{array}\right.\]
\end{Lemma}
\noindent{\bf Proof.}~~Multiplying $(\ref{1.5})_{1}$ by $p(1)-p(v)$, $(\ref{1.5})_{2}$ by $u$, and then adding the resulting equations, we have
\begin{eqnarray}\label{2.5}
\left(\Phi(v)+\frac{u^2}{2}\right)_{t}+\frac{2\nu u^2}{v^{2}}=H_{x}+\varepsilon^{2}u\left(-\frac{v_{xx}}{v^{4}}+\frac{2v_{x}^{2}}{v^{5}}\right)_{x},
\end{eqnarray}
where $H=\displaystyle(p(1)-p(v))u+2\nu\frac{u_{x}u}{v^{2}}$.

Using equation $(\ref{1.5})_{1}$, we get
\begin{eqnarray}\label{2.6}
\varepsilon^{2}u\left(-\frac{v_{xx}}{v^{4}}+\frac{2v_{x}^{2}}{v^{5}}\right)_{x}
&=&\varepsilon^{2}\left\{u\left(-\frac{v_{xx}}{v^{4}}+\frac{2v_{x}^{2}}{v^{5}}\right)\right\}_{x}
-\varepsilon^{2}u_{x}\left(-\frac{v_{xx}}{v^{4}}+\frac{2v_{x}^{2}}{v^{5}}\right)\nonumber\\
&=&\{\cdots\}_{x}-\varepsilon^{2}v_{t}\left(-\frac{v_{xx}}{v^{4}}+\frac{2v_{x}^{2}}{v^{5}}\right)\\
&=&\{\cdots\}_{x}-\varepsilon^{2}\left(\frac{v_x^2}{2v^4}\right)_{t}.\nonumber
\end{eqnarray}
Here and hereafter,  $\{\cdots\}_{x}$ denotes some terms which will disappear after integration.

Combining (\ref{2.5}) and (\ref{2.6}), and integrating the resulting equation over $[0,t]\times\mathbb{R}$, we obtain (\ref{2.4}). This completes the proof of Lemma 2.1.

From the elementary  energy estimates (\ref{2.4}), we can get the following:
\begin{Lemma}If  the conditions   of  Proposition 2.2 hold, then for any  interval $[i,i+1]\subset\mathbb{R}$ with $i\in\mathbb{Z}$, there exist $a_i(t), b_i(t)\in[i,i+1]$, and a positive constant
$C_{6}$ depending only on $\gamma,\underline{V}, \overline{V}$, and $\|(v_{0}-1, u_{0}, \varepsilon v_{0x})\|$, such that  for $t\in[0,T]$,
$$v(t,a_i(t))\leq C_6, \quad\frac{1}{v(t,b_i(t))}\leq C_6,\quad \int_{i}^{i+1}v(t,x)dx\leq C_{6}.$$
\end{Lemma}
\noindent{\bf Proof.}~~For any $i\in\mathbb{Z}$, we can obtain  from (2.4) that for $\gamma>1$,
\begin{equation*}
\int_{i}^{i+1}\left(v-1+\displaystyle\frac{1}{\gamma-1}\left(v^{-\gamma+1}-1\right)\right)dx
=\int_{i}^{i+1}\Phi(v)dx\leq C(\underline{V}, \overline{V})\|(v_0-1,u_0, \varepsilon v_{0x})\|^2,
\end{equation*}
which implies that
\begin{equation}\label{2.7}
\int_{i}^{i+1}\left(v+\displaystyle\frac{1}{\gamma-1}v^{-\gamma+1}\right)dx
\leq C(\underline{V}, \overline{V})\|(v_0-1,u_0, \varepsilon v_{0x})\|^2+\frac{\gamma}{\gamma-1}.
\end{equation}
Then if $\gamma>1$, Lemma 2.2 follows from (\ref{2.7}) and the mean value theorem  immediately.

If  $\gamma=1$,  we also have from (2.4) that for any  $i\in\mathbb{Z}$,
\begin{equation}\label{2.8-1}
\int_{i}^{i+1}\left(v-1-\ln v\right)dx
=\int_{i}^{i+1}\Phi(v)dx\leq C(\underline{V}, \overline{V})\|(v_0-1,u_0, \varepsilon v_{0x})\|^2.
\end{equation}
Since the function $\varphi(v)=v-1-\ln v$ is convex, we obtain from (\ref{2.8-1}) and Jensen's inequality:
 \[\varphi\left(\frac{\int_a^bf(x)dx}{b-a}\right)\leq\frac{\int_a^b\varphi(f(x))dx}{b-a}, \quad \mbox{for $f(x)\in L^1[a,b]$} \]
that
 \begin{equation}\label{2.8-2}
\int_{i}^{i+1}v (t,x)dx-\ln \int_{i}^{i+1}v(t,x)dx-1\leq C(\underline{V}, \overline{V})\|(v_0-1,u_0, \varepsilon v_{0x})\|^2:=C_\ast.
\end{equation}

Denoting the two positive solutions of the equation $x-\ln x-1=C_\ast$ by $\alpha$ and $\beta$, then we have from (\ref{2.8-2}) that
\begin{equation}\label{2.8}
\alpha\leq\int_{i}^{i+1}v (t,x)dx\leq \beta.
\end{equation}
 Then for the case $\gamma=1$, Lemma 2.2 can be deduced  from  (\ref{2.8}) and the mean value theorem.  This finishes the proof of Lemma 2.2.

Based on Lemma 2.2, we now deduce the lower bound of the specific volume $v(t,x)$ by using the idea developed in \cite{Jiang-Xin-Zhang}, see also \cite{Guo-Jiu-Xin,Li-Li-Xin}.
\begin{Lemma}Suppose that the conditions  of  Proposition 2.2 hold,  then there exists a positive constant
$C_{7}$ depending only on $\varepsilon, \gamma, \underline{V}, \overline{V}$, and $\|(v_{0}-1, u_{0}, \varepsilon v_{0x})\|$ such that
\begin{eqnarray}\label{2.9}
v(t,x)\geq C_{7}, \quad\forall(t,x)\in[0,T]\times\mathbb{R}.
\end{eqnarray}
\end{Lemma}
\noindent{\bf Proof.}~~For any $x\in\mathbb{R}$, there exists an integer $ i_0\in \mathbb{Z}$ such that $x\in[i_0,i_0+1]$. Then we have from Lemma 2.1 and the H\"{o}lder inequality that
\begin{eqnarray}
\frac{1}{v(t,x)}-\frac{1}{v(t,b_{i_{0}}(t))}&=&\int_{b_{i_{0}}(t)}^{x}\left(\frac{1}{v(t,y)}\right)_{y}dy\nonumber\\
&\leq&\int_{i_{0}}^{i_{0}+1}\left|\frac{v_{x}(t,x)}{v^{2}(t,x)}\right|dx\nonumber\\
&\leq&\left(\int_{i_{0}}^{i_{0}+1}\frac{v_{x}^{2}(t,x)}{v^{4}(t,x)}dx\right)^{\frac{1}{2}}\nonumber\\
&\leq&C(\underline{V}, \overline{V})\varepsilon^{-1}\|(v_{0}-1, u_{0}, \varepsilon v_{0x})\|,\nonumber
\end{eqnarray}
which, together with Lemma 2.2 gives (\ref{2.9}). The proof of Lemma 2.3 is completed.

Next, we derive  the estimate on $\displaystyle\int_{\mathbb{R}}\frac{v_{x}^{2}}{v^{4}}dx$.
\begin{Lemma}Assume that  the conditions   of  Proposition 2.2 hold,  then there exists a positive constant
$C_{8}(\underline{V}, \overline{V})$ such that
\begin{eqnarray}\label{2.10}
&&\nu\int_{\mathbb{R}}\frac{v_{x}^{2}}{v^{4}}dx
+\int_0^t\int_{\mathbb{R}}\frac{\gamma v_{x}^{2}}{v^{\gamma+3}}dxd\tau
+\varepsilon^{2}\int_0^t\int_{\mathbb{R}}\left\{\frac{1}{v^{2}}\left[\left(\frac{v_{x}}{v^{2}}\right)_{x}\right]^{2}
+\frac{v_{xx}^{2}}{v^{6}}+\frac{v_{x}^{4}}{v^{8}}\right\}dxd\tau\nonumber\\
&\leq&C_{8}(\underline{V}, \overline{V})\left(\nu\|v_{0x}\|^{2}+\|(v_{0x},u_{0})\|^{2}
+\frac{1}{\nu}\|(v_{0}-1, u_{0}, \varepsilon v_{0x})\|^{2}\right)
\end{eqnarray}
holds for all $t\in[0,T]$.
\end{Lemma}
\noindent{\bf Proof.}~~The equation   $(\ref{1.5})_{2}$ can be rewritten  as follows:
\begin{eqnarray*}
\left(\frac{2\nu v_{x}}{v^{2}}-u\right)_{t}=p'(v)v_{x}+\varepsilon^{2}\left(\frac{v_{xx}}{v^{4}}-\frac{2v_{x}^{2}}{v^{5}}\right)_{x}.
\end{eqnarray*}
Multiplying the above equation by $\displaystyle\frac{v_{x}}{v^{2}}$ yields
\begin{eqnarray}\label{2.11}
\left(\frac{\nu v_{x}^{2}}{v^{4}}-u\frac{v_{x}}{v^{2}}\right)_{t}-p'(v)\frac{v_{x}^{2}}{v^{2}}+\left(\frac{uu_x}{v^2}\right)_x=
\frac{u_{x}^{2}}{v^{2}}+\varepsilon^{2}\left(\frac{v_{xx}}{v^{4}}-\frac{2v_{x}^{2}}{v^{5}}\right)_{x}\frac{v_{x}}{v^{2}}.
\end{eqnarray}
Using direct computations, we have
\begin{eqnarray}\label{2.12}
\varepsilon^{2}\left(\frac{v_{xx}}{v^{4}}-\frac{2v_{x}^{2}}{v^{5}}\right)_{x}\frac{v_{x}}{v^{2}}&=&
\{\cdots\}_{x}-\varepsilon^{2}\left(\frac{v_{xx}}{v^{4}}-\frac{2v_{x}^{2}}{v^{5}}\right)\left(\frac{v_{x}}{v^{2}}\right)_{x}\nonumber\\
&=&\{\cdots\}_{x}-\varepsilon^{2}\left\{\frac{v_{xx}^{2}}{v^{6}}-4\frac{v_{xx}v_{x}^{2}}{v^{7}}+4\frac{v_{x}^{4}}{v^{8}}\right\}\nonumber\\
&=&\{\cdots\}_{x}-\varepsilon^{2}\left(\frac{v_{xx}}{v^{3}}-\frac{2v_{x}^{2}}{v^{4}}\right)^{2}\\
&=&\{\cdots\}_{x}-\frac{\varepsilon^{2}}{v^{2}}\left[\left(\frac{v_{x}}{v^{2}}\right)_{x}\right]^{2}.\nonumber
\end{eqnarray}
On the other hand,  it also holds that
\begin{eqnarray}\label{2.13}
\varepsilon^{2}\left(\frac{v_{xx}}{v^{4}}-\frac{2v_{x}^{2}}{v^{5}}\right)_{x}\frac{v_{x}}{v^{2}}
&=&\{\cdots\}_{x}-\varepsilon^{2}\frac{v_{xx}^{2}}{v^{6}}+4\varepsilon^{2}\frac{\left(\frac{1}{3}v_{x}^{3}\right)_{x}}{v^{7}}-4\varepsilon^{2}\frac{v_{x}^{4}}{v^{8}}\\
&=&\{\cdots\}_{x}-\frac{\varepsilon^{2}v_{xx}^{2}}{v^{6}}+\frac{16\varepsilon^{2}}{3}\frac{v_{x}^{4}}{v^{8}}.\nonumber
\end{eqnarray}

Substituting (\ref{2.12}) into (\ref{2.11}),  and integrating the resultant equation over $[0,t]\times\mathbb{R}$ gives rise to
\begin{eqnarray}\label{2.14}
&&\frac{\nu}{2}\int_{\mathbb{R}}\frac{v_{x}^{2}}{v^{4}}dx+\int_0^t\int_{\mathbb{R}}\frac{\gamma v_{x}^{2}}{v^{\gamma+3}}dxd\tau
+\varepsilon^{2}\int_0^t\int_{\mathbb{R}}\frac{1}{v^{2}}\left[\left(\frac{v_{x}}{v^{2}}\right)_{x}\right]^{2}dxd\tau\nonumber\\
&\leq&\frac{\nu}{2}\int_{\mathbb{R}}\frac{v_{0x}^{2}}{v_{0}^{4}}dx+\frac{1}{2}\int_{\mathbb{R}}\frac{v_{0x}^{2}}{v_{0}^{4}}dx
+\frac{1}{2}\int_{\mathbb{R}}u_{0}^{2}dx+\frac{1}{2\nu}\int_{\mathbb{R}}u^{2}dx+\int_0^t\int_{\mathbb{R}}\frac{u_{x}^{2}}{v^{2}}dxd\tau,
\end{eqnarray}
where the following Cauchy inequality:
\begin{eqnarray*}
\int_{\mathbb{R}}\frac{v_{x}}{v^{2}}u\,dx\leq\frac{\nu}{2}\int_{\mathbb{R}}\frac{v_{x}^{2}}{v^{4}}dx+\frac{1}{2\nu}\int_{\mathbb{R}}u^{2}dx
\end{eqnarray*}
has been used. Then it follows from (\ref{2.14}) and Lemma 2.1 that
\begin{eqnarray}\label{2.15}
&&\frac{\nu}{2}\int_{\mathbb{R}}\frac{v_{x}^{2}}{v^{4}}dx
+\int_0^t\int_{\mathbb{R}}\frac{\gamma v_{x}^{2}}{v^{\gamma+3}}dxd\tau+\varepsilon^{2}\int_0^t\int_{\mathbb{R}}\frac{1}{v^{2}}\left[\left(\frac{v_{x}}{v^{2}}\right)_{x}\right]^{2}dxd\tau
\nonumber\\
&\leq&C(\underline{V}, \overline{V})\left(\nu\|v_{0x}\|^{2}+\|(v_{0x},u_{0})\|^{2}
+\frac{1}{\nu}\|(v_{0}-1, u_{0}, \varepsilon v_{0x})\|^{2}\right).
\end{eqnarray}

Similarly, inserting  (\ref{2.13}) into (\ref{2.11}), and integrating the resultant equation over $[0,t]\times\mathbb{R}$, we obtain from Germain and LeFloch's inequality (\ref{5.1}), the Cauchy inequality, and Lemma 2.1 that
\begin{eqnarray}\label{2.16}
&&\frac{\nu}{2}\int_{\mathbb{R}}\frac{v_{x}^{2}}{v^{4}}dx
+\int_0^t\int_{\mathbb{R}}\frac{\gamma v_{x}^{2}}{v^{\gamma+3}}dxd\tau
+\frac{\varepsilon^{2}}{49}\int_0^t\int_{\mathbb{R}}\frac{v_{xx}^{2}}{v^{6}}dxd\tau
+\frac{\varepsilon^{2}}{9}\int_0^t\int_{\mathbb{R}}\frac{v_{x}^{4}}{v^{8}}dxd\tau\nonumber\\
&\leq&C(\underline{V}, \overline{V})\left(\nu\|v_{0x}\|^{2}+\|(v_{0x},u_{0})\|^{2}
+\frac{1}{\nu}\|(v_{0}-1, u_{0}, \varepsilon v_{0x})\|^{2}\right).
\end{eqnarray}
Then (\ref{2.10}) follows from (\ref{2.15}) and (\ref{2.16}) immediately. This ends the proof of Lemma 2.4.

\begin{Remark}
From Lemmas 2.1 and 2.4, we can get  the following Bresch-Desjardins type entropy inequality developed by Bresch and Desjardins \cite{Bresch-Desjardins,Bresch-Desjardins-Lin} for the multi-dimensional viscous compressible flow:
\begin{eqnarray*}
&&\int_{\mathbb{R}}\left(\Phi(v)+\frac{1}{2}\left(u-\frac{2\nu v_x}{v^2}\right)^2+\frac{\varepsilon^{2}}{2}\frac{v_x^2}{v^{4}}\right)dx
+\int_0^t\int_{\mathbb{R}}\frac{\gamma v_{x}^{2}}{v^{\gamma+3}}dxd\tau\nonumber\\
&&+\varepsilon^{2}\int_0^t\int_{\mathbb{R}}\left\{\frac{1}{v^{2}}\left[\left(\frac{v_{x}}{v^{2}}\right)_{x}\right]^{2}
+\frac{v_{xx}^{2}}{v^{6}}+\frac{v_{x}^{4}}{v^{8}}\right\}dxd\tau\nonumber\\
&\leq&C(\underline{V}, \overline{V},\varepsilon, \nu,\|v_{0}-1\|, \|u_{0}\|),\quad \forall \,t\in[0,T].
\end{eqnarray*}
Such type inequality was first obtained by Kanel' \cite{Kanel} for the one-dimensional case.
\end{Remark}

Based on Lemmas 2.1 and 2.4, we can get the following key lemma, which will paly an important role in establishing the upper bound of the specific volume $v(t,x)$.
\begin{Lemma}If the conditions  of Proposition 2.2 hold, then there exists a positive constant
$C_{9}(T)$ depending only on $\nu,\varepsilon,\gamma,\underline{V}, \overline{V}$, $\|(v_{0}-1, u_{0}, \varepsilon v_{0x})\|$, and $T$,  such that for  $t\in[0,T]$,
\begin{equation}\label{2.17}
\int_0^t\left\|\frac{u(\tau)}{v^{\frac{3}{4}}(\tau)}\right\|^2_{L^\infty}d\tau\leq C_{9}(T).
\end{equation}
\end{Lemma}
\noindent{\bf Proof.}~~First,  we have  from the H\"{o}lder inequality that
\begin{eqnarray}\label{2.18}
\frac{u^2(t,x)}{v^{\frac{3}{2}}(t,x)}&=&2\int_{-\infty}^x\frac{u(t,y)}{v^{\frac{3}{4}}(t,y)} \left(\frac{u(t,y)}{v^{\frac{3}{4}}(t,y)}\right)_ydy\nonumber\\
&\leq&C\int_{\mathbb{R}}\left(\left|\frac{u(t,x) u_x(t,x)}{v^{\frac{3}{2}}(t,x)}\right|+\left|\frac{u^2(t,x) v_x(t,x)}{v^{\frac{5}{2}}(t,x)}\right|\right)dx\nonumber\\
&\leq&C\left(\|u(t)\|\left\|\frac{u_x(t)}{v^{\frac{3}{2}}(t)}\right\| +\|u(t)\|^2\left\|\frac{v_x(t)}{v^{\frac{5}{2}}(t)}\right\|_{L^\infty}\right).
\end{eqnarray}
Similarly, it  holds that
\begin{eqnarray}\label{2.19}
\frac{v^2_x(t,x)}{v^{5}(t,x)}&=&\int_{-\infty}^x2\frac{v_y(t,y)}{v^{\frac{5}{2}}(t,y)} \left(\frac{v_y(t,y)}{v^{\frac{5}{2}}(t,y)}\right)_ydy\nonumber\\
&\leq&C\int_{\mathbb{R}}\left(\left|\frac{v_x(t,x)v_{xx}(t,x)}{v^{5}(t,x)}\right| +\left|\frac{v^3_x(t,x)}{v^{6}(t,x)}\right|\right)dx\nonumber\\
&\leq&C\left(\left\|\frac{v_{x}(t)}{v^{2}(t)}\right\| \left\|\frac{v_{xx}(t)}{v^{3}(t)}\right\|
+\left\|\frac{v_{x}(t)}{v^{2}(t)}\right\|\left\|\frac{v^2_{x}(t)}{v^{4}(t)}\right\|\right),
\end{eqnarray}
which implies that
\begin{equation}\label{2.20}
\left\|\frac{v_x(t)}{v^{\frac{5}{2}}(t)}\right\|_{L^\infty}\leq C\left\|\frac{v_{x}(t)}{v^{2}(t)}\right\|^{\frac{1}{2}}\left\|\left(\frac{v_{xx}(t)}{v^{3}(t)}, \frac{v^2_{x}(t)}{v^{4}(t)}\right)\right\|^{\frac{1}{2}}.
\end{equation}

Combining (\ref{2.18}) and (\ref{2.20}), then we have from  Lemmas 2.1 and 2.3,  and the Young inequality that
\begin{equation}\label{2.21}
\begin{split}
\left\|\frac{u(t)}{v^{\frac{3}{4}}(t)}\right\|^2_{L^\infty}&\leq C\left(\|u(t)\|\left\|\frac{u_x(t)}{v(t)}\right\|\left\|\frac{1}{v}\right\|^{\frac{1}{2}}_{L^\infty_{T,x}} +\|u(t)\|^2\left\|\frac{v_{x}(t)}{v^{2}(t)}\right\|^{\frac{1}{2}}
\left\|\left(\frac{v_{xx}(t)}{v^{3}(t)}, \frac{v^2_{x}(t)}{v^{4}(t)}\right)\right\|^{\frac{1}{2}}\right)\\
&\leq C_{10}\left(\left\|\frac{u_x(t)}{v(t)}\right\|^2+\left\|\left(\frac{v_{xx}(t)}{v^{3}(t)}, \frac{v^2_{x}(t)}{v^{4}(t)}\right)\right\|^{2}+1\right),
\end{split}
\end{equation}
where $C_{10}$ is a positive constant depending only on $\varepsilon, \gamma, \underline{V}, \overline{V}$, and $\|(v_{0}-1, u_{0}, \varepsilon v_{0x})\|$.

Integrating (\ref{2.21}) with respect to $t$ over $(0,t)$,  and then employing Lemmas 2.1 and 2.4,  we have (\ref{2.17}) holds. The proof of Lemma 2.5 is finished.

In order to deduce the upper bound of the specific volume $v(t,x)$, we assume that  $\nu\geq\varepsilon$ and let $c_{\pm}:=\nu\pm\sqrt{\nu^{2}-\varepsilon^2}$ be two positive constants. Introducing a new effective velocity $\xi(t,x)$ by
 \begin{equation*}
 \xi(t,x):=u(t,x)-c_{+}\frac{v_{x}(t,x)}{v^2(t,x)},
 \end{equation*}
then the unknown functions $(v(t,x), \xi(t,x))$ satisfy  the following Cauchy problem:
\begin{eqnarray}\label{2.22}
\left\{\begin{array}{ll}
     v_{t}-\xi_{x}=\displaystyle c_{+}\left(\frac{v_{x}}{v^{2}}\right)_{x},\quad x\in\mathbb{R},\,\,t>0,\\[4mm]
      \xi_{t}+\displaystyle\frac{\gamma}{c_{+}}v^{1-\gamma}\xi=\displaystyle\frac{\gamma}{c_{+}}v^{1-\gamma}u+c_{-}\left(\frac{\xi_{x}}{v^{2}}\right)_{x}
\end{array}\right.
\end{eqnarray}
with initial data
\begin{eqnarray}\label{2.23}
(v(t,x),\xi(t,x))|_{t=0}=\left(v_{0}(x),\xi_{0}(x)\right)=\displaystyle
\left(v_{0}(x),u_{0}(x)-c_{+}\frac{v_{0x}(x)}{v^{2}_{0}(x)}\right), \,\, x\in\mathbb{R}.
\end{eqnarray}

The derivation of the system (\ref{2.22}) can be found in the Appendix. The following lemma gives the $L^\infty$ norm on the effective  velocity $\xi(t,x)$.
\begin{Lemma}Suppose that the conditions of  Proposition 2.2 hold, then  there exists a positive constant
$C_{11}(T)$ depending only on $\nu,\varepsilon,\gamma,\underline{V}, \overline{V}$, $\|v_{0}-1\|_2$, $\|u_{0}\|_1$, and $T$, such that for  $t\in[0,T]$,
\begin{equation}\label{2.24}
\|\xi(t)\|_{L^{\infty}}\leq C_{11}(T).
\end{equation}
\end{Lemma}
\noindent{\bf Proof.}~~Taking $(\ref{2.22})_{2}\times\xi|\xi|^{p}$ with $p\geq0$, and integrating the resulting  equation with respect to $t$ and $x$ over $[0,t]\times\mathbb{R}$  gives
\begin{eqnarray}\label{2.25}
&&\frac{1}{p+2}\int_{\mathbb{R}}|\xi|^{p+2}dx+\frac{\gamma}{c_{+}}\int_{0}^{t}\int_{\mathbb{R}}v^{1-\gamma}|\xi|^{p+2}dxd\tau
+c_{-}(p+1)\int_{0}^{t}\int_{\mathbb{R}}\frac{|\xi|^{p}\xi_{x}^{2}}{v^{2}}dxd\tau\nonumber\\
&=&\frac{1}{p+2}\int_{\mathbb{R}}|\xi_{0}|^{p+2}dx+\underbrace{\frac{\gamma}{c_{+}}\int_{0}^{t}\int_{\mathbb{R}}v^{1-\gamma}u\xi|\xi|^{p}dxd\tau}_{\mathcal{R}_0}.
\end{eqnarray}
Now we estimate the reminder term $\mathcal{R}_0$. Since
\begin{eqnarray*}
\left\|v^{1-\gamma}(t)u(t)\right\|_{L^{p+2}}^{p+2}&=&\int_{\mathbb{R}}(v(t,x))^{(1-\gamma)(p+2)} \left\|\frac{u(t)}{v^{\frac{3}{4}}(t)}\right\|_{L^{\infty}}^{p}v^{\frac{3}{4}p}(t,x)u^{2}(t,x)dx\\
&\leq&\left\|\frac{1}{v}\right\|_{L^{\infty}_{T,x}}^{(p+2)\left(\gamma-1-\frac{3p}{4(p+2)}\right)} \left\|\frac{u(t)}{v^{\frac{3}{4}}(t)}\right\|_{L^{\infty}}^{p}\|u(t)\|^{2},
\end{eqnarray*}
it follows from Lemmas 2.1 and 2.3  that
\begin{eqnarray}\label{2.26}
\left\|v^{1-\gamma}(t)u(t)\right\|_{L^{p+2}}&\leq& \left\|\frac{1}{v}\right\|_{L^{\infty}_{T,x}}^{(\gamma-1)} \left\|\frac{u(t)}{v^{\frac{3}{4}}(t)}\right\|_{L^{\infty}}^{\frac{p}{p+2}}
\|u(t)\|^{\frac{2}{p+2}}\nonumber\\
&\leq&C_{12}\left(\varepsilon,\gamma, \underline{V}, \overline{V}, \|(v_{0}-1, u_{0}, \varepsilon v_{0x})\|\right)\left\|\frac{u(t)}{v^{\frac{3}{4}}(t)}\right\|_{L^{\infty}}^{\frac{p}{p+2}},
\end{eqnarray}
where we have used the assumption that $\gamma\geq\frac{7}{4}$ such that the exponent $\gamma-1-\frac{3p}{4(p+2)}\geq0$.

Utilizing (\ref{2.26}), the H\"older inequality, the Young inequality, and Lemma 2.5, we have
\begin{eqnarray}\label{2.27}
\mathcal{R}_0&\leq&\frac{\gamma}{c_{+}}\int_{0}^{t}\|\xi(\tau)\|_{L^{p+2}}^{p+1} \left\|v^{1-\gamma}(\tau)u(\tau)\right\|_{L^{p+2}}d\tau\nonumber\\
&\leq&\frac{\gamma C_{12}}{c_{+}}\int_{0}^{t}\left(1+\|\xi(\tau)\|_{L^{p+2}}^{p+2}\right) \left(\left\|\frac{u(\tau)}{v^{\frac{3}{4}}(\tau)}\right\|_{L^{\infty}}^{2}+1\right)d\tau\\
&\leq&\frac{\gamma C_{12}}{c_{+}}C(T)+\frac{\gamma C_{12}}{c_{+}}\int_{0}^{t}\left(\left\|\frac{u(\tau)}{v^{\frac{3}{4}}(\tau)}\right\|_{L^{\infty}}^{2}+1\right)
\|\xi(\tau)\|_{L^{p+2}}^{p+2}d\tau\nonumber.
\end{eqnarray}

Inserting (\ref{2.27}) into (\ref{2.25}), then we have from the Gronwall inequality that
\begin{eqnarray*}
\|\xi(t)\|_{L^{p+2}}^{p+2}\leq
\left(\|\xi_{0}\|_{L^{p+2}}^{p+2}+\frac{\gamma C_{12}}{c_{+}}C(T)(p+2)\right)e^{\frac{\gamma C_{12}}{c_{+}}(p+2)(T+C_{9}(T))},
\end{eqnarray*}
which implies that
\begin{eqnarray}\label{2.28}
\|\xi(t)\|_{L^{p+2}}\leq
2^{\frac{1}{p+2}}\left(\|\xi_{0}\|_{L^{p+2}}
+\left(\frac{\gamma C_{12}}{c_{+}}C(T)(p+2)\right)^{\frac{1}{p+2}}\right)e^{\frac{\gamma C_{12}}{c_{+}}(T+C_{9}(T))}.
\end{eqnarray}

Letting $p\rightarrow+\infty$, then (\ref{2.24}) follows from (\ref{2.28}) and the fact that $\|\xi_0\|_{L^\infty}\leq C(\|v_0-1\|_2+\|u_0\|_1)$ immediately. This completes the proof of Lemma 2.6.

Now we drive the upper bound of the specific volume $v(t,x)$ in the following lemma.
\begin{Lemma}If the conditions   of  Proposition 2.2 hold, then  there exists a positive constant
$C_{13}(T)$ depending only on $\nu,\varepsilon,\gamma,\underline{V}, \overline{V}$, $\|v_{0}-1\|_2$, $\|u_{0}\|_1$, and $T$, such that
\begin{equation}\label{2.29}
v(t,x)\leq  C_{13}(T), \,\,\,\,\forall\,(t,x)\in[0,T]\times\mathbb{R}.
\end{equation}
\end{Lemma}
\noindent{\bf Proof.}~~Letting  $\phi(t,x)=v(t,x)-1$ and multiplying $(\ref{2.22})_{1}$ by $\phi|\phi|^{2}$, then integrating the resultant equation over $[0,t]\times\mathbb{R}$, we have
\begin{eqnarray}\label{2.30}
\frac{1}{4}\int_{\mathbb{R}}|\phi|^{4}dx+3c_{+}\int_{0}^{t}\int_{\mathbb{R}}|\phi|^{2}\frac{\phi_{x}^{2}}{v^{2}}dxd\tau
=\frac{1}{4}\int_{\mathbb{R}}|\phi_{0}|^{4}dx+\underbrace{\int_{0}^{t}\int_{\mathbb{R}}\phi|\phi|^{2}\xi_{x}dxd\tau}_{\mathcal{R}_1},
\end{eqnarray}
where $\phi_0(x)=v_0(x)-1$.

Using integration by parts, the Cauchy inequality, and the Young equality, we obtain
\begin{eqnarray}\label{2.31}
\mathcal{R}_1&=&-3\int_{0}^{t}\int_{\mathbb{R}}\phi^{2}\phi_{x}\xi dxd\tau\nonumber\\
&\leq&\frac{3c_{+}}{2}\int_{0}^{t}\int_{\mathbb{R}}\phi^{2}\frac{\phi_{x}^{2}}{v^{2}} dxd\tau+\frac{3}{2c_{+}}\int_{0}^{t}\int_{\mathbb{R}}\phi^{2}v^{2}\xi^{2} dxd\tau\\
&\leq&\frac{3c_{+}}{2}\int_{0}^{t}\int_{\mathbb{R}}\phi^{2}\frac{\phi_{x}^{2}}{v^{2}} dxd\tau+\frac{3}{2c_{+}}\int_{0}^{t}\int_{\mathbb{R}}2\phi^{2}(\phi^{2}+1)\xi^{2} dxd\tau\nonumber\\
&\leq&\frac{3c_{+}}{2}\int_{0}^{t}\int_{\mathbb{R}}\phi^{2}\frac{\phi_{x}^{2}}{v^{2}} dxd\tau+\frac{3}{c_{+}}\int_{0}^{t}\int_{\mathbb{R}}\left(\phi^{4}\xi^{2}+\frac{1}{2}\phi^{4}+\frac{1}{2}\xi^{4}\right) dxd\tau\nonumber\\
&\leq&\frac{3c_{+}}{2}\int_{0}^{t}\int_{\mathbb{R}}\phi^{2}\frac{\phi_{x}^{2}}{v^{2}}dxd\tau
+\frac{3}{c_{+}}\int_{0}^{t}\|\phi(\tau)\|_{L^{4}}^{4}\|\xi(\tau)\|_{L^{\infty}}^{2}d\tau\nonumber\\
&&+\frac{3}{2c_{+}}\int_{0}^{t}\left(\|\phi(\tau)\|_{L^{4}}^{4}+\|\xi(\tau)\|^{2}\|\xi(\tau)\|_{L^{\infty}}^{2}\right)d\tau\nonumber\\
&\leq&\frac{3c_{+}}{2}\int_{0}^{t}\int_{\mathbb{R}}\phi^{2}\frac{\phi_{x}^{2}}{v^{2}}dxd\tau
+\frac{9(C_{11}(T))^{2}}{2c_{+}}\int_{0}^{t}\|\phi(\tau)\|_{L^{4}}^{4}d\tau
+\frac{3}{2c_{+}}C_{14}(T),\nonumber
\end{eqnarray}
where $C_{14}(T)$ is a positive constant depending only on $\nu,\varepsilon,\gamma,\underline{V}, \overline{V}$, $\|v_{0}-1\|_2$, $\|u_{0}\|_1$, and $T$, and we have used Lemma 2.6 and the fact that
\begin{eqnarray*}
\|\xi(t)\|^{2}&\leq&2\left(\|u(t)\|^{2}+c_{+}^{2}\left\|\frac{v_{x}(t)}{v^{2}(t)}\right\|^{2}\right)\\
&\leq&C\left(\|(v_{0}-1, u_{0}, \varepsilon v_{0x})\|, \nu,\varepsilon, \underline{V}, \overline{V}\right)
\end{eqnarray*}
due to Lemma 2.1.

Substituting (\ref{2.31}) into (\ref{2.30}), and  making use of the Gronwall inequality, we have
\begin{eqnarray*}
\|\phi(t)\|_{L^{4}}^{4}&\leq&\left(\|\phi_{0}\|_{L^{4}}^{4}+\frac{6}{c_{+}}C_{14}(T)\right)e^{\frac{18(C_{11}(T))^{2}}{c_{+}}T},
\end{eqnarray*}
which leads to that for any $i\in\mathbb{Z}$,
\begin{equation}\label{2.32}
\int_{i}^{i+1}v^{4}dx\leq16\int_{i}^{i+1}(\phi^{4}+1)dx \leq C_{15}(T),
\end{equation}
where $C_{15}(T)$ is a positive constant depending only on $\nu,\varepsilon,\gamma,\underline{V}, \overline{V}$, $\|v_{0}-1\|_2$, $\|u_{0}\|_1$, and $T$.

Now we prove the upper bound of $v(t,x)$ by using the idea developed in \cite{Jiang-Xin-Zhang}, see also \cite{Guo-Jiu-Xin,Li-Li-Xin}. Indeed, for any $x\in\mathbb{R}$, there exists an integer $i_0\in \mathbb{Z}$ such that $x\in[i_{0},i_{0}+1]$. Then we obtain from  the H\"older inequality, (\ref{2.4}), and (\ref{2.32}) that
\begin{eqnarray*}
v(t,x)-v(t,a_{i_0}(t))&=&\int_{a_{i_0}(t)}^{x}v_{y}(t,y)dy\\
&\leq&\left(\int_{i_{0}}^{i_{0}+1}\frac{v_{x}(t,x)^{2}}{v(t,x)^{4}}dx\right)^{\frac{1}{2}}
\left(\int_{i_{0}}^{i_{0}+1}v(t,x)^{4}dx\right)^{\frac{1}{2}}\\
&\leq&\varepsilon^{-1}C(\underline{V}, \overline{V})\|(v_{0}-1, u_{0}, \varepsilon v_{0x})\|(C_{15}(T))^{\frac{1}{2}},
\end{eqnarray*}
which together with Lemma 2.2 gives (\ref{2.29}). This finishes the proof of Lemma 2.7.

As a consequence of Lemmas 2.1, 2.3, 2.4, and 2.7, we get the following:
\begin{Corollary}Assume that  the conditions  of Proposition 2.2 hold, then there exists a positive constant
$C_{16}(T)$ depending only on $\nu,\varepsilon,\gamma,\underline{V}, \overline{V}$, $\|v_{0}-1\|_2$, $\|u_{0}\|_1$, and $T$, such that for  $t\in[0,T]$,
\begin{equation}\label{2.33}
\|v(t)-1\|_1^{2}+\|u(t)\|^{2}+\int_{0}^{t}\left(\|v_{x}(\tau)\|_1^{2}+\|u_{x}(\tau)\|^{2}\right)d\tau
\leq C_{16}(T)\left(\|v_{0}-1\|_{1}^{2}+\|u_{0}\|^{2}\right).
\end{equation}
\end{Corollary}

The next lemma gives the estimate on $\|u_x(t)\|$.
\begin{Lemma}If the conditions  of  Proposition 2.2 hold, then there exists a positive constant
$C_{17}(T)$ depending only on $\nu,\varepsilon,\gamma,\underline{V}, \overline{V}$, $\|v_{0}-1\|_2$, $\|u_{0}\|_1$, and $T$, such that for  $t\in[0,T]$,
\begin{equation}\label{2.34}
\|(u_{x}(t), v_{xx}(t))\|^{2}+\int_{0}^{t}\|u_{xx}(\tau)\|^{2}d\tau\leq C_{17}(T)(\|v_{0}-1\|_{2}^{2}+\|u_{0}\|_{1}^{2}).
\end{equation}
\end{Lemma}
\noindent{\bf Proof.}~~Multiplying $(\ref{1.5})_{2}$ by $-u_{xx}$, and using $(\ref{1.5})_{1}$, we obtain
\begin{eqnarray*}
&&\left(\frac{u_{x}^{2}}{2}+\frac{\varepsilon^{2}v_{xx}^{2}}{2v^{4}}\right)_{t}+\frac{2\nu u_{xx}^{2}}{v^{2}}\\
&=&u_{xx}v_{x}p'(v)+\frac{4\nu u_{x}v_{x}u_{xx}}{v^{3}}-\frac{2\varepsilon^{2}v^2_{xx}u_{x}}{v^{5}}
-2\varepsilon^{2}\left(\frac{2v_{x}v_{xx}}{v^{5}}-\frac{5v_{x}^{3}}{v^{6}}\right)u_{xx}+\{\cdots\}_{x}.
\end{eqnarray*}
Integrating the above equation over $[0,t]\times\mathbb{R}$, we get from  Lemmas 2.3 and 2.7 that
\begin{eqnarray}\label{2.35}
\|(u_{x}(t),\varepsilon v_{xx}(t))\|^{2}+\int_{0}^{t}\|u_{xx}(\tau)\|^{2}d\tau
\leq C(T)\left(\|u_{0x}\|^{2}+\varepsilon^{2}\|v_{0xx}\|^{2}+K_{1}+K_{2}\right),
\end{eqnarray}
where
\begin{eqnarray*}
&&K_{1}=\int_{0}^{t}\int_{\mathbb{R}}\varepsilon^{2}|v_{xx}^{2}u_{x}|dxd\tau,\\
&&K_{2}=\int_{0}^{t}\int_{\mathbb{R}}\left(|v_{x}|+|u_{x}v_{x}|+\varepsilon^{2}|v_{x}v_{xx}|+\varepsilon^{2}|v_{x}^{3}|\right)|u_{xx}|dxd\tau.
\end{eqnarray*}
From the Sobolev inequality, the Cauchy inequality, and Corollary 2.1, we have 
\begin{eqnarray}\label{2.36}
K_{1}&\leq&\varepsilon^{2}\int_{0}^{t}\|u_{x}(\tau)\|^{\frac{1}{2}}\|u_{xx}(\tau)\|^{\frac{1}{2}}\|v_{xx}(\tau)\|^{2}d\tau\nonumber\\
&\leq&\eta\int_{0}^{t}\|u_{xx}(\tau)\|^{2}d\tau+C_{\eta}\int_{0}^{t}\|u_{x}(\tau)\|^{\frac{2}{3}}\|v_{xx}(\tau)\|^{\frac{8}{3}}d\tau\\
&\leq&\eta\int_{0}^{t}\|u_{xx}(\tau)\|^{2}d\tau+C_{\eta}\int_{0}^{t}\|v_{xx}(\tau)\|^{4}d\tau +C_{\eta}\int_{0}^{t}\|u_{x}(\tau)\|^{2}d\tau\nonumber
\end{eqnarray}
and
\begin{eqnarray}\label{2.37}
K_{2}&\leq&\eta\int_{0}^{t}\|u_{xx}(\tau)\|^{2}d\tau+C_{\eta}\int_{0}^{t}\left(\|v_{x}(\tau)\|^{2} +\|u_{x}(\tau)\|\|u_{xx}(\tau)\|\|v_{x}(\tau)\|^{2}\right.\nonumber\\
&&\left.+\|v_{x}(\tau)\|\|v_{xx}(\tau)\|^{3}+\|v_{x}(\tau)\|^{4}\|v_{xx}(\tau)\|^{2}\right)d\tau\nonumber\\
&\leq&\eta\int_{0}^{t}\|u_{xx}(\tau)\|^{2}d\tau+C_{\eta}\int_{0}^{t}\left\{\|v_{x}(\tau)\|^{2} +\|u_{x}(\tau)\|\|u_{xx}(\tau)\|\sup_{0\leq\tau\leq t}\{\|v_{x}(\tau)\|^{2}\}\right.\nonumber\\
&&\left.+\|v_{x}(\tau)\|_{1}^{2}\|v_{xx}(\tau)\|^{2}
+\sup_{0\leq\tau\leq t}\{\|v_{x}(\tau)\|^{2}\}\|v_{x}(\tau)\|^{2}\|v_{xx}(\tau)\|^{2}\right\}d\tau\\
&\leq&2\eta\int_{0}^{t}\|u_{xx}(\tau)\|^{2}d\tau+C_{\eta}(T)\int_{0}^{t}\left(\|(u_{x}(\tau),v_{x}(\tau))\|^{2}
+\|v_{x}(\tau)\|_{1}^{2}\cdot\|v_{xx}(\tau)\|^{2}\right)d\tau\nonumber.
\end{eqnarray}
Here and hereafter, $\eta>0$ is a small positive constant, and $C_\eta$ is a constant depending on $\eta$.

Then (\ref{2.34}) follows from (\ref{2.35})-(\ref{2.37}), the smallness of $\eta$,  Corollary 2.1 and the Gronwall inequality. The proof of Lemma 2.8 is completed.

Finally, we estimate $\|v_{xx}(t)\|$.
\begin{Lemma} Let the conditions   of  Proposition 2.2 hold, then  there exists a positive constant
$C_{18}(T)$ which depends only on $\nu,\varepsilon,\gamma,\underline{V}, \overline{V}$, $\|v_{0}-1\|_2$, $\|u_{0}\|_1$, and $T$, such that
\begin{equation}\label{2.38}
\|v_{xx}(t)\|^{2}+\int_{0}^{t}\|v_{xx}(\tau)\|_{1}^{2}d\tau\leq C_{18}(T)(\|v_{0}-1\|_{2}^{2}+\|u_{0}\|_{1}^{2})
\end{equation}
holds for  all $t\in[0,T]$.
\end{Lemma}
\noindent{\bf Proof.}~~Multiplying $(\ref{1.5})_{2}$ by $v_{xxx}$ and using equation $(\ref{1.5})_{1}$,  we have
\begin{eqnarray}\label{2.39}
&&\left(\frac{\nu v_{xx}^{2}}{v^{2}}-v_{xx}u_{x}\right)_{t}+\frac{\varepsilon^{2}v_{xxx}^{2}}{v^{4}}-p'(v)v_{xx}^{2}\nonumber\\
&=&u_{xx}^{2}+\left(p''(v)v_{x}^{2}v_{xx}+\frac{2\nu v_{xx}^{2}u_{x}}{v^{3}}+\frac{8\nu v_{x}v_{xx}u_{xx}}{v^{3}}-\frac{12\nu v_{x}^{2}v_{xx}u_{x}}{v^{4}}\right)\nonumber\\
&&+\varepsilon^{2}\left(\frac{8v_{x}v_{xx}}{v^{5}}-\frac{10v_{x}^{3}}{v^{6}}\right)v_{xxx}+\{\cdots\}_{x}\nonumber\\
&:=&u_{xx}^{2}+K_{3}+K_{4}+\{\cdots\}_{x}.
\end{eqnarray}

Integrating (\ref{2.39}) with respect to $t$ and $x$ over $[0,t]\times\mathbb{R}$,  we derive from Lemmas 2.3 and 2.7 that
\begin{eqnarray}\label{2.40}
&&\nu\|v_{xx}(t)\|^{2}+\int_{0}^{t}\|v_{xx}(\tau)\|^{2}d\tau+\varepsilon^{2}\int_{0}^{t}\|v_{xxx}(\tau)\|^{2}d\tau\nonumber\\
&\leq&C(T)\left(\nu\|v_{0xx}\|^{2}+\int_{0}^{t}\|u_{xx}(\tau)\|^{2}d\tau+\left|\int_{\mathbb{R}}v_{xx}u_{x}dx\right|+
\left|\int_{\mathbb{R}}v_{0xx}u_{0x}dx\right|\right.\nonumber\\
&&\left.+\int_{0}^{t}\int_{\mathbb{R}}(|K_{3}|+|K_{4}|)dxd\tau\right).
\end{eqnarray}
The Cauchy inequality and the Sobolev inequality imply that
\begin{eqnarray}\label{2.41}
\left|\int_{\mathbb{R}}v_{xx}u_{x}dx\right|\leq \eta\|v_{xx}(t)\|^{2}+C_{\eta}\|u_{x}(t)\|^{2},
\end{eqnarray}
\begin{eqnarray}\label{2.42}
\int_{0}^{t}\int_{\mathbb{R}}|K_{3}|dxd\tau&\leq&C(T)\int_{0}^{t}\int_{\mathbb{R}}
\left(|v_{x}^{2}v_{xx}|+|v_{xx}^{2}u_{x}|+|v_{x}v_{xx}u_{xx}|+|v_{x}^{2}v_{xx}u_{x}|\right)dxd\tau\nonumber\\
&\leq&C(T)\int_{0}^{t}\|v_{x}(\tau)\|_{L^{\infty}}^{2}\|(v_{x}(\tau), u_{x}(\tau), u_{xx}(\tau))\|\|v_{xx}(\tau)\|d\tau\nonumber\\
&&+C(T)\int_{0}^{t}\|u_{x}(\tau)\|^{\frac{1}{2}}\|u_{xx}(\tau)\|^{\frac{1}{2}}\|v_{xx}(\tau)\|^{2}d\tau\\
&\leq&C(T)\int_{0}^{t}\left(\|(v_{x}(\tau),v_{xx}(\tau),u_{x},u_{xx}(\tau))\|^{2}+\|v_{xx}(\tau)\|^{4}\right)d\tau\nonumber\\
&\leq&C(T)\int_{0}^{t}\|(v_{x}(\tau),u_{x}(\tau))\|_{1}^{2}d\tau+C(T)\int_{0}^{t}\|v_{xx}(\tau)\|^{4}d\tau,\nonumber
\end{eqnarray}
and
\begin{eqnarray}\label{2.43}
\int_{0}^{t}\int_{\mathbb{R}}|K_{4}|dxd\tau&\leq&\frac{\varepsilon^{2}}{2}\int_{0}^{t}\|v_{xxx}(\tau)\|^{2}d\tau
+C(T)\varepsilon^{2}\int_{0}^{t}\int_{\mathbb{R}}\left(|v_{x}^{2}v_{xx}^{2}|+v_{x}^{6}\right)dxd\tau\nonumber\\
&\leq&\frac{\varepsilon^{2}}{2}\int_{0}^{t}\|v_{xxx}(\tau)\|^{2}d\tau
+C(T)\int_{0}^{t}\left(\|v_{x}(\tau)\|_{L^\infty}^{2}\|v_{xx}(\tau)\|^2+
\|v_{x}(\tau)\|_{L^\infty}^{4}\|v_{x}(\tau)\|^2\right)d\tau\nonumber\\
&\leq&\frac{\varepsilon^{2}}{2}\int_{0}^{t}\|v_{xxx}(\tau)\|^{2}d\tau+C(T)\int_{0}^{t}\|v_{x}(\tau)\|_{1}^{2}d\tau.
\end{eqnarray}
Then (\ref{2.38}) can be obtained by combining (\ref{2.40})-(\ref{2.43}), the smallness of $\eta$, Corollary 2.1, and Lemma 2.8. This ends the proof of Lemma 2.9.

Combining Corollary 2.1 and Lemmas 2.8-2.9,  we have
\begin{Corollary}Suppose that  the conditions   of  Proposition 2.2 hold, then  there exists a positive constant
$C_{19}(T)$ depending only on $\nu,\varepsilon,\gamma,\underline{V}, \overline{V}$, $\|v_{0}-1\|_2$, $\|u_{0}\|_1$, and $T$,  such that for  $t\in[0,T]$,
\begin{equation}\label{2.44}
\|v(t)-1\|_2^{2}+\|u(t)\|_1^{2}+\int_{0}^{t}\left(\|v_{x}(\tau)\|_2^{2}+\|u_{x}(\tau)\|_1^{2}\right)d\tau
\leq C_{19}(T)\left(\|v_{0}-1\|_{1}^{2}+\|u_{0}\|^{2}\right).
\end{equation}
\end{Corollary}

Then by  the standard continuation argument, one can extend the local solution of the Cauchy problem (\ref{1.5})-(\ref{1.6}) step by step to be a global one. Furthermore, by employing a large-time stability analysis,  we can obtain the uniform-in-time positive lower and upper bound of the specific volume $v(t,x)$ in the following lemma.
\begin{Lemma} Let the conditions of Proposition 2.2 hold, then there exists a positive constant $C_{0}$ which depends  only on $\nu,\varepsilon, \gamma, \underline{V}, \overline{V}, \|v_0-1\|_2$, and $\|u_0\|_1$, such that
 \begin{equation}\label{2.45}
 C_{0}^{-1}\leq v(t,x)\leq C_{0},\quad\forall\,(t,x)\in[0,\infty)\times\mathbb{R}.
 \end{equation}
\end{Lemma}\noindent{\bf Proof.}~~ Since the proof is similar to that of Lemma 4.3 in \cite{Chen-Li-2021},  we only give the sketch here. To prove $(\ref{2.45})$,  one only needs to show that
\begin{equation}\label{2.46}
    \lim_{t\rightarrow+\infty}\left\|\frac{1}{v(t)}-1\right\|_{L^{\infty}}=0.
\end{equation}
In fact, suppose  (\ref{2.46}) holds for the temporary, then there exists a sufficient large constant $T_0>0$ such that
  \[
    \left|\frac{1}{v(t,x)}-1\right| \leq \left\|\frac{1}{v(t)}-1\right\|_{L^{\infty}} < \frac{1}{2},\quad\forall\, (t,x)\in[T_0,+\infty)\times\mathbb{R},
  \]
which implies that
  \[
    \frac{2}{3}\leq v(t,x)\leq 2,\quad\forall\, (t,x)\in[T_0,+\infty)\times\mathbb{R}.
  \]
On the other hand,  we obtain from Lemmas 2.3 and 2.7 that
 \[C_{7}\leq v(t,x)\leq C_{13}(T_0),\quad\forall\,(t,x)\in[0,T_0]\times\mathbb{R}.\]
 Then (\ref{2.45}) follows immediately by setting  $C_{0}=\max\left\{C_{13}(T_0),C_7^{-1}, 2\right\}$.

Then it remains to prove (\ref{2.46}) holds. First, notice that the estimates in Lemmas 2.1 and 2.4 are uniform-in-time, and  the global-in-time  solutions for the Cauchy problem (\ref{1.5})-(\ref{1.6}) has been obtained. Therefore, we can get from Lemmas 2.1 and 2.4 that
\begin{eqnarray}\label{2.47}
&&\sup_{0\leq t<+\infty}\int_{\mathbb{R}}\left(\Phi(v)+\frac{u^2}{2}+\frac{\varepsilon^2v_x^2}{2v^{4}}\right)dx+
      2\nu\int_0^{+\infty}\int_{\mathbb{R}}\frac{u_x^2}{v^{2}}dxd\tau\nonumber\\[2mm]
      &\leq&C_{20}(\underline{V}, \overline{V})\|(v_0-1,u_0,\varepsilon v_{0x})\|^2
 \end{eqnarray}
and
\begin{eqnarray}\label{2.48}
&&\sup_{0\leq t<+\infty}\int_{\mathbb{R}}\frac{\nu v_{x}^{2}}{v^{4}}dx
+\int_0^{+\infty}\int_{\mathbb{R}}\frac{\gamma v_{x}^{2}}{v^{\gamma+3}}dxd\tau
+\varepsilon^{2}\int_0^{+\infty}\int_{\mathbb{R}}\left\{\frac{1}{v^{2}}\left[\left(\frac{v_{x}}{v^{2}}\right)_{x}\right]^{2}
+\frac{v_{xx}^{2}}{v^{6}}+\frac{v_{x}^{4}}{v^{8}}\right\}dxd\tau\nonumber\\
&\leq&C_{8}(\underline{V}, \overline{V})\left(\nu\|v_{0x}\|^{2}+\|(v_{0x},u_{0})\|^{2}
+\frac{1}{\nu}\|(v_{0}-1, u_{0}, \varepsilon v_{0x})\|^{2}\right).
\end{eqnarray}

Based on (\ref{2.47})-(\ref{2.48}), one can show (\ref{2.46}) holds  by using the same method as that of Lemma 4.3 in \cite{Chen-Li-2021}, the details are omitted  here for brevity.  This finishes  the proof of Lemma 2.10.\\

With the above preparations in hand, we now turn to prove Proposition 2.2. \\
\noindent{\bf Proof of Proposition 2.2.}~~Since the uniform-in-time lower and upper bounds of the specific volume $v(t,x)$ has been obtained in Lemma 2.10, one can get from (\ref{2.47})-(\ref{2.48}) that
\begin{equation}\label{2.49}
\|v(t)-1\|_1^{2}+\|u(t)\|^{2}+\int_{0}^{t}\left(\|v_{x}(\tau)\|_1^{2}+\|u_{x}(\tau)\|^{2}\right)d\tau
\leq C_{21}\left(\|v_{0}-1\|_{1}^{2}+\|u_{0}\|^{2}\right),\,\,\forall\, t\in[0,T],
\end{equation}
where $C_{21}$ is a positive constant depending only on $\nu,\varepsilon, \gamma, \underline{V}, \overline{V}, \|v_0-1\|_2$, and $\|u_0\|_1$.

Furthermore, similar to the proof of Lemmas 2.8 and 2.9,  we  can  also obtain for any $t\in[0,T]$ that
\begin{equation}\label{2.50}
\left\|\left(u_{x}(t), v_{xx}(t)\right)\right\|^{2}+\int_{0}^{t}\|u_{xx}(\tau)\|^{2}d\tau\leq C_{22}(\|v_{0}-1\|_{2}^{2}+\|u_{0}\|_{1}^{2})
\end{equation}
and
\begin{equation}\label{2.51}
\|v_{xx}(t)\|^{2}+\int_{0}^{t}\|v_{xx}(\tau)\|_{1}^{2}d\tau\leq C_{23}(\|v_{0}-1\|_{2}^{2}+\|u_{0}\|_{1}^{2}),
\end{equation}
where $C_{22}, C_{23}$ are positive constants depending only on $\nu,\varepsilon, \gamma, \underline{V}, \overline{V}, \|v_0-1\|_2$, and $\|u_0\|_1$.

Then (\ref{2.2}) and (\ref{2.3}) follows from Lemma 2.10 and (\ref{2.49})-(\ref{2.51}) immediately. This completes the proof of Proposition 2.2.

\subsection{Proof of Theorem 1.1}

We are now in a position to finish  the proof of Theorem 1.1 as follows. Indeed, based on Propositions 2.1 and 2.2, one can get the global existence of strong solutions to the Cauchy problem (\ref{1.5})-(\ref{1.6}) by the standard continuation argument. Then the uniform-in-time estimates (\ref{1.7}) and (\ref{1.8}) follow directly from (\ref{2.2}) and (\ref{2.3}), respectively. Next, we turn to prove the large-time behavior (\ref{1.9}).  To this end, we infer from  (\ref{1.8}) and the original system (\ref{1.5}) that
\begin{equation}\label{2.52}
\int_0^{+\infty}\left(\|(v_x(t), u_x(t), v_{xx}(t))\|^2+\left|\frac{d}{dt}\|(v_x(t), u_x(t), v_{xx}(t))\|^2\right|\right)\,dt<\infty,
\end{equation}
which implies the following time-asymptotic behavior of solutions:
\begin{equation}\label{2.53}
\|(v_x(t), u_x(t), v_{xx}(t))\|\rightarrow0,\quad\mbox{as}\,\, t\rightarrow+\infty.
\end{equation}
Furthermore, it follows from (\ref{2.53}), (\ref{1.8}), and the Sobolev inequality that
\begin{eqnarray}\label{2.54}
 \|(v(t)-1, u(t))\|_{L^\infty}
&\leq&\|(v(t)-1,u(t))\|^{\frac{1}{2}}\|(v_x(t), u_x(t))\|^{\frac{1}{2}}\nonumber\\
&\leq&C \|(v_x(t),u_x(t))\|^{\frac{1}{2}}\rightarrow0,\quad\mbox{as}\,\, t\rightarrow+\infty.
\end{eqnarray}
Then (\ref{2.53}) together with (\ref{2.54}) leads to (\ref{1.9}) immediately.  Hence the proof of Theorem 1.1 is completed.

\section{Global  existence in the parabolic regime}
\setcounter{equation}{0}
This section is devoted to proving Theorem 1.2, which is a  global existence result in the parabolic regime $\varepsilon\ll\nu$ with large initial data.  First of all, similar to Proposition 2.1, the local existence result for the initial data $(v_0(x)-1,u_0(x))\in H^{5}(\mathbb{R})\times H^4(\mathbb{R})$  can be  stated as follows.

\begin{Proposition} [Local existence]
Suppose that  the conditions  of  Theorem 1.2 hold, and  the initial data $(v_0(x)-1,u_0(x))\in H^{5}(\mathbb{R})\times H^4(\mathbb{R})$, then there exists a sufficiently small positive constant $\tilde{t}_0=\tilde{t}_0(\underline{V},\tilde{N}_0)$ depending only on $\nu, \gamma,\underline{V}$, and $\tilde{N}_0:=(\|v_0-1\|^2_{5}+\|u_0\|^2_4)^{\frac{1}{2}}$, such that the Cauchy  problem (\ref{1.5})-(\ref{1.6}) admits a unique  solution $(v(t,x), u(t,x))\in X_4(0,\tilde{t}_0; \frac{1}{2}\underline{V}, 2\overline{V})$ and
 \[\sup_{[0,\tilde{t}_0]}\left\{\|v(t)-1\|^2_{5}+\|u(t)\|^2_{4}+\varepsilon^2\|v_x(t)\|^2_{4}\right\}+
\int_0^{\tilde{t}_0}\left(\|(u_x(\tau),v_x(\tau))\|^2_{4}+\varepsilon^2\|v_{xx}(\tau)\|^2_{4}\right)d\tau\leq \tilde{b}\tilde{N}_0^2,
\]
 where $\tilde{b}>1$ is a positive constant depending only on $\underline{V}, \overline{V}$.
\end{Proposition}

Then in order to prove Theorem 1.2, it suffices to derive the following:
\begin{Proposition} [{\it A priori} estimates]
 Assume that the conditions  of Theorem 1.2 hold,  and $(v(t,x), u\\(t,x))\in X_4(0,T;m,M)$ is a solution of the Cauchy problem  (\ref{1.5})-(\ref{1.6}) defined in  $\Pi_T=[0,T]\times\mathbb{R}$ for some positive constants $m, M$, and $T$. Moreover, $(v(t,x),u(t,x))$ satisfies  the following {\it a priori} assumption:
 \begin{eqnarray}\label{3.1}
\sup_{[0,T]}\left\{\|v(t)-1\|^2_{5}+\|u(t)\|^2_{4}+\varepsilon^2\|v_x(t)\|^2_{4}\right\}+
\int_0^T\left(\|(u_x(\tau),v_x(\tau))\|^2_{4}+\varepsilon^2\|v_{xx}(\tau)\|^2_{4}\right)d\tau\leq N^2
\end{eqnarray}
 for some positive constant $N>1$. Then  if $\varepsilon$ is  sufficiently small such that $\varepsilon^2M^{\gamma-1}N^4<1$, it  holds that
\begin{equation}\label{3.2}
C_{24}^{-1}\leq v(t,x), \quad \forall\,(t,x)\in[0,T]\times\mathbb{R},
\end{equation}
\begin{equation}\label{3.3}
\sup_{\mathbb{R}}v(t,x)\leq \left(\sup_{\mathbb{R}}v(s,x)+C_{25}\right)e^{C(\nu,\gamma)(t-s)},\quad\forall\,0\leq s\leq t\leq T,
\end{equation}
\begin{eqnarray}\label{3.4}
&&\|v(t)-1\|^2_{5}+\|u(t)\|^2_{4}+\varepsilon^2\|v_x(t)\|^2_{4}+
\int_0^t\left(\|(u_x(\tau),v_x(\tau))\|^2_{4}+\varepsilon^2\|v_{xx}(\tau)\|^2_{4}\right)d\tau\nonumber\\
&&\leq C_{26}(M)\left(\|v_0-1\|^2_{5}+\|u_0\|^2_{4}\right), \quad \forall\, t\in[0,T],
\end{eqnarray}
where the constant $C(\nu,\gamma)=\max\{\frac{\gamma}{2\nu(2-\gamma)}, \frac{1}{\nu}\}$, $C_{24}$ and $C_{25}$ are two  positive constants depending only on $\nu, \gamma,\underline{V}, \overline{V}, \|v_0-1\|_{1}$, and $\|u_0\|$, and $C_{26}(M)$ is a positive constant depending only on $\nu, \gamma,\underline{V}, \overline{V}, \|v_0-1\|_{5}$, $\|u_0\|_4$, and $M$.
\end{Proposition}

\subsection{{\it A priori} estimates}
This subsection is devoted to proving Proposition 3.2, which  can be obtained by a series of lemmas below. First, notice that Lemmas 2.1, 2.2 and 2.4 still hold for the case when $\varepsilon$ is sufficiently small, and the constant $C_6$ given in Lemma  2.2 can depend only on  $\gamma,\underline{V}, \overline{V}$, and $\|(v_{0}-1, u_{0},  v_{0x})\|$ if $\varepsilon\ll1$.
Since our goal here is to show the uniform-in-$\varepsilon$ estimates (\ref{3.2})-(\ref{3.4}), the lower bound estimate of the specific volume $v(t,x)$ stated in Lemma 2.3 are not valid any longer. Instead of Lemma 2.3, we establish the following lemma.
\begin{Lemma}Let the conditions   of  Proposition 3.2 hold, then  there exists a positive constant
$C_{24}$ depending only on $\gamma, \nu,\underline{V}, \overline{V}$, and $\|(v_{0}-1, u_{0},  v_{0x})\|$ such that
\begin{eqnarray}\label{3.5}
v(t,x)\geq C_{24}^{-1}, \quad\forall(t,x)\in[0,T]\times\mathbb{R}.
\end{eqnarray}
\end{Lemma}
\noindent{\bf Proof.}~~First, we have from (\ref{2.10}) with $\varepsilon\ll1$ that
\begin{eqnarray}\label{3.6}
&&\nu\int_{\mathbb{R}}\frac{v_{x}^{2}}{v^{4}}dx
+\int_0^t\int_{\mathbb{R}}\frac{\gamma v_{x}^{2}}{v^{\gamma+3}}dxd\tau
+\varepsilon^{2}\int_0^t\int_{\mathbb{R}}\left\{\frac{1}{v^{2}}\left[\left(\frac{v_{x}}{v^{2}}\right)_{x}\right]^{2}
+\frac{v_{xx}^{2}}{v^{6}}+\frac{v_{x}^{4}}{v^{8}}\right\}dxd\tau\nonumber\\
&\leq&C_{8}(\underline{V}, \overline{V})\left(\nu\|v_{0x}\|^{2}+\|(v_{0x},u_{0})\|^{2}
+\frac{1}{\nu}\|(v_{0}-1, u_{0}, v_{0x})\|^{2}\right):=C_{27},\quad \forall t\in[0,T].
\end{eqnarray}
Then for any $x\in\mathbb{R}$, there exists an integer $ i_0\in \mathbb{Z}$ such that $x\in[i_0,i_0+1]$. It follows from  (\ref{3.6}) and the H\"{o}lder inequality that
\begin{eqnarray}
\frac{1}{v(t,x)}-\frac{1}{v(t,b_{i_{0}}(t))}&=&\int_{b_{i_{0}}(t)}^{x}\left(\frac{1}{v(t,y)}\right)_{y}dy
\leq\int_{i_{0}}^{i_{0}+1}\left|\frac{v_{x}(t,x)}{v^{2}(t,x)}\right|dx\nonumber\\
&\leq&\left(\int_{i_{0}}^{i_{0}+1}\frac{v^{2}_{x}(t,x)}{v^{4}(t,x)}dx\right)^{\frac{1}{2}}\nonumber\\
&\leq&\left(\frac{C_{27}}{\nu}\right)^{\frac{1}{2}},\nonumber
\end{eqnarray}
which together with Lemma 2.2 implies  (\ref{3.5}). This ends the proof of Lemma 3.1.

For the upper bound of the specific volume $v(t,x)$,  we have the  following lemma.
\begin{Lemma}
Under the conditions   of  Proposition 3.2,  there exists a positive constant
$C_{25} $ depending only on $\nu, \gamma,\underline{V}, \overline{V}, \|v_0-1\|_{1}$, and $\|u_0\|$, such that if   $\varepsilon^2M^{\gamma-1}N^4<1$, then
(\ref{3.3}) holds for all $0\leq s\leq t\leq T$.
\end{Lemma}
\noindent{\bf Proof.} In order to derive the upper bound of  $v(t,x)$ for the parabolic regime,  we introduce a new effective velocity $\omega(t,x)$ by
\begin{equation}\label{3.7}
\omega(t,x)=u(t,x)-2\nu\frac{v_{x}(t,x)}{v^{2}(t,x)},
\end{equation}
then the Cauchy problem (\ref{1.5})-(\ref{1.6}) can be rewritten as
\begin{eqnarray}\label{3.8}
\left\{\begin{array}{ll}
    \displaystyle v_{t}-\omega_{x}=\displaystyle2\nu\left(\frac{v_{x}}{v^{2}}\right)_{x}, \quad x\in\mathbb{R},\,\,t>0,\\[2mm]
    \omega_{t}+\displaystyle\frac{\gamma}{2\nu}v^{1-\gamma}\omega=\displaystyle\frac{\gamma}{2\nu}v^{1-\gamma}u
    +\varepsilon^{2}\left(-\frac{v_{xx}}{v^{4}}+\frac{2v_{x}^{2}}{v^{5}}\right)_{x}
\end{array}\right.
\end{eqnarray}
with the initial data
\begin{eqnarray}\label{3.9}
(v(t,x),\omega(t,x))|_{t=0}=(v_{0}(x),\omega_{0}(x))=\left(v_{0}(x), u_{0}(x)-2\nu\frac{v_{0x}(x)}{v^2_{0}(x)}\right), \quad x\in \mathbb{R}.
\end{eqnarray}

 Now we  show the upper bound estimate of $v(t,x)$ by modifying  the argument from the proof of Theorem 2.1  in \cite{Burtea-Haspot}. Indeed, differentiating $(\ref{3.8})_2$ with respect to $x$ once,  and setting $\tilde{\omega}=\omega_x$,  then we obtain
 \begin{equation}\label{3.10}
   \tilde{\omega}_t+\frac{\gamma}{2\nu} v^{1-\gamma}\tilde{\omega}=\gamma(1-\gamma)\frac{v_x^2}{v^{\gamma+2}}
   +\frac{\gamma}{2\nu} v^{1-\gamma}u_x+\varepsilon^{2}\left(-\frac{v_{xx}}{v^{4}}+\frac{2v_{x}^{2}}{v^{5}}\right)_{xx}.
  \end{equation}

Let $\tilde{p}=\tilde{\omega}+F(v)$  with $F(v)$ being  a function to be determined later. As  \cite{Burtea-Haspot}, $\tilde{p}$ is called {\it the effective pressure}.  Then one can obtain from (\ref{3.10}) and $(\ref{1.5})_1$ that
 \begin{eqnarray}\label{3.11}
   \tilde{p}_t
   =\left(F^\prime(v)+\frac{\gamma}{2\nu} v^{1-\gamma}\right)u_x-\frac{\gamma}{2\nu} v^{1-\gamma}\tilde{\omega}+\gamma(1-\gamma)\frac{v_x^2}{v^{2+\gamma}}+\varepsilon^{2}\left(-\frac{v_{xx}}{v^{4}}+\frac{2v_{x}^{2}}{v^{5}}\right)_{xx}.
 \end{eqnarray}
Define the function $F(v)$ by
 \begin{equation}\label{3.12}
   F(v)=\left\{\begin{array}{ll}
   \displaystyle\frac{\gamma}{2\nu(\gamma-2)}v^{2-\gamma}, \quad if\, 1\leq\gamma<2, \\[4mm]
  \displaystyle-\frac{\gamma}{2\nu}\ln v, \quad \quad if\, \gamma=2,
 \end{array}\right.
 \end{equation}
 then it follows  from (\ref{3.11}) that for $1\leq\gamma<2$,
 \begin{equation}\label{3.13}
   \tilde{p}_t+\frac{\gamma}{2\nu} v^{1-\gamma}\tilde{p}=\frac{\gamma^2}{4\nu^2(\gamma-2)}v^{3-2\gamma}
 +\gamma(1-\gamma)\frac{v_x^2}{v^{\gamma+2}}+\varepsilon^{2}\left(-\frac{v_{xx}}{v^{4}}+\frac{2v_{x}^{2}}{v^{5}}\right)_{xx}.
 \end{equation}

Since $\tilde{p}(t,x)$ is bounded and continuous in $[0,T]\times\mathbb{R}$,  and  $\lim_{x\rightarrow\pm\infty}\tilde{p}(t,x)=F(1)$,  we have for each $t\in[0,T]$  that
\begin{equation}\label{3.14}
\tilde{p}_M(t):=\sup_{\mathbb{R}}\tilde{p}(t,x)\geq F(1).
\end{equation}
The function $\tilde{p}_M(t)$ is continuous in $[0,T]$ due to the fact that $\tilde{p}(t,x)\in C([0,T]\times\mathbb{R})$. Therefore,  the set
\[D:=\left\{t\geq0|\tilde{p}_M(t)>F(1)\right\}\]
is open in $[0,T]$. Consequently, it holds that
$
D=I_0\bigcup(\bigcup_{j\in\mathbb{Z}^+}I_j),
$
where $I_j=(a_j, b_j)$ ($j\in\mathbb{Z}^+$) are disjoint open intervals, and $I_0=\emptyset$ if  $\tilde{p}_M(0)=F(1)$, while $I_0=[0,a_0)$ for some number $a_0\in[0,T]$ if  $\tilde{p}_M(0)>F(1)$. Moreover, from the definition of $I_j(j\in\mathbb{Z}^+)$,  we have that $\tilde{p}_M(a_j)=\tilde{p}_M(b_j)=F(1)$.

For any $t\notin D$, we have
\begin{eqnarray}\label{3.15}
\sup_{\mathbb{R}}\tilde{p}(t,x)=\tilde{p}_M(t)\leq F(1)=\frac{\gamma}{2\nu(\gamma-2)}\leq0,\quad if \,1\leq\gamma<2.
\end{eqnarray}

For any $t\in D$, we now show that $\tilde{p}_M(t)$  has  also a upper bound.  Indeed, if $t\in D$, then there exists $j_0\in\mathbb{N}$ such that $t\in I_{j_0}$.  Since $\tilde{p}(t,x)$ is continuous with respect to $x$ in $\mathbb{R}$, and $\tilde{p}_M(t):=\sup_{\mathbb{R}}\tilde{p}(t,x)>F(1)$, and $\lim_{x\rightarrow\pm\infty} \tilde{p}(t,x)=F(1)$, there exists a point $x_t\in\mathbb{R}$ such that
\[ \tilde{p}_M(t):=\sup_{\mathbb{R}}\tilde{p}(t,x)=\tilde{p}(t,x_t).\]

We claim that $\tilde{p}_M(t)$ is differentiable in $I_{j_0}$ and $(\tilde{p}_M)^\prime(t)=\partial_t\tilde{p}(t,x_t)$.

In fact, for $t\in I_{j_0}$ and $t+\Delta t\in I_{j_0}$, we have
\begin{eqnarray*}
\tilde{p}_M(t+\Delta t)-\tilde{p}_M(t)&=&\sup_{\mathbb{R}}\tilde{p}(t+\Delta t, x)-\sup_{\mathbb{R}}\tilde{p}(t, x)\nonumber\\
&=&\sup_{\mathbb{R}}\left(\tilde{p}(t+\Delta t, x)-\tilde{p}(t, x)\right)\nonumber\\
&=&\sup_{\mathbb{R}}\tilde{p}_t(t+\theta\Delta t, x)\Delta t\nonumber
\end{eqnarray*}
with $0<\theta<1$ being a constant. Therefore, it holds that
\[\left(\tilde{p}_M\right)^\prime(t)=\lim_{\Delta t\rightarrow0}\frac{\tilde{p}_M(t+\Delta t)-\tilde{p}_M(t)}{\Delta t}=\sup_{\mathbb{R}}\tilde{p}_t(t, x).\]
From equation (\ref{3.13}) and the a priori assumptions $m\leq v(t,x)\leq M$ and (\ref{3.1}),  we get that  $\|\tilde{p}_t(t)\|_{L^\infty}\leq C(\nu, \varepsilon,\gamma, m, M, N)<+\infty$. Therefore,  $\sup_{\mathbb{R}}\tilde{p}_t(t, x)$ exists and $\tilde{p}_M(t)$ is differentiable in $I_{j_0}$.

Next, we prove that $(\tilde{p}_M)^\prime(t)=\partial_t\tilde{p}(t,x_t)$.  On one hand,
\begin{eqnarray}\label{3.16}
   (\tilde{p}_M)^\prime(t)&=&\lim_{h\rightarrow0^+}\frac{ \tilde{p}_M(t+h)-\tilde{p}_M(t)}{h}=\lim_{h\rightarrow0^+}\frac{ \tilde{p}(t+h,x_{t+h})-\tilde{p}(t,x_t)}{h}\nonumber\\
   &\geq&\lim_{h\rightarrow0^+}\frac{ \tilde{p}(t+h,x_{t})-\tilde{p}(t,x_t)}{h}=\partial_t\tilde{p}(t,x_t).
 \end{eqnarray}
On the other hand,
\begin{eqnarray}\label{3.17}
   (\tilde{p}_M)^\prime(t)&=&\lim_{h\rightarrow0^+}\frac{ \tilde{p}_M(t)-\tilde{p}_M(t-h)}{h}=\lim_{h\rightarrow0^+}\frac{ \tilde{p}(t,x_{t})-\tilde{p}(t-h,x_{t-h})}{h}\nonumber\\
   &\leq&\lim_{h\rightarrow0^+}\frac{ \tilde{p}(t,x_{t})-\tilde{p}(t-h,x_t)}{h}=\partial_t\tilde{p}(t,x_t).
\end{eqnarray}
Here in (\ref{3.16}) and (\ref{3.17}), we have used the inequalities $\tilde{p}(t+h,x_{t+h})\geq\tilde{p}(t+h,x_{t})$ and $\tilde{p}(t-h,x_{t-h})\geq\tilde{p}(t-h,x_{t})$, respectively.

Then combining  (\ref{3.16}) and (\ref{3.17}) gives that  $(\tilde{p}_M)^\prime(t)=\partial_t\tilde{p}(t,x_t)$. Thus the claim is proved.

Letting $x=x_t$ in (\ref{3.13}), then we have
  \begin{equation}\label{3.18}
 (\tilde{p}_M)^\prime(t)+\frac{\gamma}{2\nu} v^{1-\gamma}(t,x_t)\tilde{p}_M(t)\leq \varepsilon^{2}\left(-\frac{v_{xx}}{v^{4}}+\frac{2v_{x}^{2}}{v^{5}}\right)_{xx}(t,x_t),
 \end{equation}
 where we have used the assumption that $1\leq\gamma<2$.

Using Lemma 3.1, the {\it a priori} assumption (\ref{3.1}) and the Sobolev inequality, we have
\begin{equation}\label{3.19}
\begin{split}
\left|\varepsilon^{2}\left(-\frac{v_{xx}}{v^{4}}+\frac{2v_{x}^{2}}{v^{5}}\right)_{xx}(t,x_t)\right|& =\varepsilon^2\left|\left(-\frac{v_{xxxx}}{v^4}
+\frac{v_xv_{xxx}}{v^5}+\frac{8v_{xx}^2}{v^5}-\frac{70v_{xx}{v_x^2}}{v^6}+\frac{60{v_x^4}}{v^7}\right)(t,x_t)\right|\\
&\leq C\varepsilon^2(C_{24})^{7}\left(\|v_{xxxx}(t)\|_1+(\|v_{x}(t)\|^2+1)\|v_{xx}(t)\|^2\right)\\
&\quad+C\varepsilon^2(C_{24})^{7}(\|v_{xx}(t)\|+\|v_{xxxx}(t)\|)\|v_{xxx}(t)\|\\
&\leq C\varepsilon^2N^4.
\end{split}
\end{equation}
Consequently, we deduce from (\ref{3.18})-(\ref{3.19}), the Gronwall inequality, and the a priori assumption $v(t,x)\leq M$  for all $(t,x)\in[0,T]\times\mathbb{R}$ that
 \begin{eqnarray}\label{3.20}
(\tilde{p}_M)(t)&\leq& \displaystyle\tilde{p}_M(a_{j_0})e^{-\frac{\gamma}{2\nu} \int_{a_{j_0}}^tv^{1-\gamma}(s,x_s)ds}+C\int_{a_{j_0}}^t \varepsilon^2N^4e^{-\frac{\gamma}{2\nu}\int_s^tv^{1-\gamma}(\tau,x_\tau)d\tau}ds\nonumber\\
&\leq& \displaystyle\tilde{p}_M(a_{j_0})+C\varepsilon^2N^4\int_{a_{j_0}}^t e^{-\frac{\gamma}{2\nu}M^{1-\gamma}(t-s)}ds\nonumber\\
&\leq&F(1)+C\frac{2\nu}{\gamma}\varepsilon^2M^{\gamma-1}N^4\leq C_{28},\quad \forall\, t\in I_{j_0},
  \end{eqnarray}
provided that $\varepsilon>0$ is sufficiently small such that $\varepsilon^2M^{\gamma-1}N^4<1$,  where $C_{28}$ is a positive constant depending only on $\nu, \gamma,\underline{V}, \overline{V}, \|v_0-1\|_1$, and $\|u_0\|$,  and we have used the fact that $\tilde{p}_M(a_{j_0})=F(1)$.

Moreover, for the case of $\gamma=2$, one can verify that (\ref{3.20}) also holds by employing the same argument as above. Then we obtain from (\ref{3.15}) and (\ref{3.20}) that for any $(t,x)\in[0,T]\times\mathbb{R}$,
\begin{equation}\label{3.21}
\omega_x(t,x)\leq C_{28}+\left\{
  \begin{array}{ll}
    \displaystyle\frac{\gamma}{2\nu(2-\gamma)}v^{2-\gamma},\quad if\,1\leq\gamma<2,\\[3mm]
   \displaystyle\frac{\gamma}{2\nu}\ln v, \quad if\,\gamma=2.
  \end{array}
  \right.
\end{equation}

Next, taking into account the fact that $\lim_{x\rightarrow\pm\infty}v(t,x)=1$ and $v(t,x)\in C([0,T]\times\mathbb{R})$,  we set $v_M(t)=\sup_{\mathbb{R}}v(t,x)$,  and consider the set
\[E=\left\{t\geq0|v_M(t)>1\right\}=A_0\bigcup(\bigcup_{j\in\mathbb{Z}^+} A_j),\]
where for $j\in\mathbb{Z}^+$, $A_j=(\tilde{a}_j,\tilde{b}_j)$ are open disjoint intervals,   and $A_0=\emptyset$ if  $v_M(0):=\sup_{\mathbb{R}}v(0,x)=1$, while $A_0=[0,\tilde{a}_0)$ for some number $\tilde{a}_0\in[0,T]$ if  $v_M(0)>1$. Following the same argument as above, we can get that in any $A_j(j\geq0)$, there exists a point $x_t\in\mathbb{R}$ such that
\[v_M(t):=\sup_{\mathbb{R}}v(t,x)=v(t,x_t).\]
Then as previously, it follows from  $(\ref{3.8})_1$ and (\ref{3.21})  that for any $t\in  E$,
\begin{equation}\label{3.22}
\begin{split}
(v_M)^\prime(t)&\leq \omega_x(t,x_t)\leq C_{28}+\left\{
  \begin{array}{ll}
     \displaystyle\frac{\gamma}{2\nu(2-\gamma)}v^{2-\gamma}(t,x_t),\quad if\,1\leq\gamma<2,\\[3mm]
   \displaystyle\frac{\gamma}{2\nu}\ln v(t,x_t), \quad if\,\gamma=2,
  \end{array} \right.\\
  &\leq C_{28}+C(\nu,\gamma)v_M(t)
  \end{split}
\end{equation}
with $C(\nu,\gamma)=\max\{\frac{\gamma}{2\nu(2-\gamma)}, \frac{1}{\nu}\}$, where we have used the fact that  $v_{xx}(t,x_t)\leq0$, $v_{x}(t,x_t)=0$ since  $v(t,x)$ attains its maximum at the point $x_t$.

Then for any $0\leq s\leq t\leq T$, we  have the following four cases:

{\it Case I:} $t,s\in A_{j_0}\subset E$ with $j_0\in\mathbb{N}$. In this case, we have from (\ref{3.22}) and the Gronwall inequality that
 \begin{equation}\label{3.23}
\sup_{\mathbb{R}}v(t,x)\leq \left(\sup_{\mathbb{R}}v(s,x)+\frac{C_{28}}{C(\nu,\gamma)}\right)e^{C(\nu,\gamma)(t-s)}.
\end{equation}

{\it Case II:} $t,s\in E$ and $t\in A_{j_1}, s\in A_{j_0}$ with $j_1>j_0$. In this case, we have
  \begin{eqnarray}\label{3.24}
\sup_{\mathbb{R}}v(t,x)&\leq& \left(\sup_{\mathbb{R}}v(\tilde{a}_{j_1},x)+\frac{C_{28}}{C(\nu,\gamma)}\right)e^{C(\nu,\gamma)(t-\tilde{a}_{j_1})}\nonumber\\
&\leq&\left(1+\frac{C_{28}}{C(\nu,\gamma)}\right)e^{C(\nu,\gamma)(t-s)},
\end{eqnarray}
where we have used the fact that $\sup_{\mathbb{R}}v(\tilde{a}_{j},x)=\sup_{\mathbb{R}}v(\tilde{b}_{j},x)=1$ for any $j\in\mathbb{N}$.

{\it Case III:} $t\in E$ and $s\notin E$. In this case, we assume  that
 $t\in A_{j_0}$ with $j_0\geq1$, then we have
  \begin{equation}\label{3.25}
  \begin{split}
\sup_{\mathbb{R}}v(t,x)&\leq \left(\sup_{\mathbb{R}}v(\tilde{a}_{j_0},x)+\frac{C_{28}}{C(\nu,\gamma)}\right)e^{C(\nu,\gamma)(t-\tilde{a}_{j_0})}\\
&\leq\left(1+\frac{C_{28}}{C(\nu,\gamma)}\right)e^{C(\nu,\gamma)(t-s)}.
\end{split}
\end{equation}

{\it Case IV:} $t\notin E$. In this case, it is easy to get
 \begin{equation}\label{3.26}
\sup_{\mathbb{R}}v(t,x)\leq1\leq \left(1+\frac{C_{28}}{C(\nu,\gamma)}\right)e^{C(\nu,\gamma)(t-s)}.
\end{equation}
Combining  (\ref{3.23})-(\ref{3.26}) and letting $C_{25}=1+\frac{C_{28}}{C(\nu,\gamma)}$, then we have  (\ref{3.3}) immediately. The proof of Lemma 3.2 is completed.

Combining  Lemmas 2.1 and 2.4, and using the a priori assumption that $
  v(t,x)\leq M$ for all $(t,x)\in[0,T]\times\mathbb{R}$, we have the following:
\begin{Corollary}
Assume that  the conditions  of  Proposition 3.2 hold,  then there exists a positive constant
$C_{29}(M)$ depending only on $\nu,\gamma,\underline{V}, \overline{V}$, $\|v_{0}-1\|_1$,  $\|u_0\|$, and $M$,  such that for $t\in[0,T]$,
\begin{eqnarray}\label{3.27}
&&\|v(t)-1\|_1^{2}+\|u(t)\|^{2}+\varepsilon^{2}\|v_{x}(t)\|^{2}+\int_{0}^{t}\left(\|(v_x(\tau),u_{x}(\tau))\|^{2}+\varepsilon^{2}\|v_{xx}(\tau)\|^{2} +\varepsilon^{2}\|v_{x}(\tau)\|_{L^{4}}^{4}\right)d\tau\nonumber\\
&\leq& C_{29}(M)\|(v_{0}-1, u_{0}, v_{0x})\|^{2}.
\end{eqnarray}
\end{Corollary}

Next, we derive the estimate on $\|u_x(t)\|$.
\begin{Lemma}
Let the conditions  of  Proposition 3.2 hold,  then  there exists a positive constant
$C_{30}(M)$ depending only on $\nu,\gamma,\underline{V}, \overline{V}$, $\|v_{0}-1\|_1$, $\|u_0\|$, and $M$, such that for $t\in[0,T]$,
\begin{eqnarray}\label{3.28}
\|(u_{x}(t),\varepsilon v_{xx}(t))\|^{2}+\int_{0}^{t}\|u_{xx}(\tau)\|^{2}d\tau\leq C_{30}(M)\left(\|v_{0}-1\|_{2}^{2}+\|u_{0}\|_{1}^{2}\right).
\end{eqnarray}
\end{Lemma}
\noindent{\bf Proof.}~~Similar to (\ref{2.35}), we have
\begin{eqnarray}\label{2.29}
\|(u_{x}(t),\varepsilon v_{xx}(t))\|^{2}+\int_{0}^{t}\|u_{xx}(\tau)\|^{2}d\tau
\leq C(M)\left(\|u_{0x}\|^{2}+\varepsilon^{2}\|v_{0xx}\|^{2}+K_{1}+K_{2}\right),
\end{eqnarray}
where $C(M)$ is a positive constant depending only on $\nu,\gamma,\underline{V}, \overline{V}$, $\|v_{0}-1\|_1$, $\|u_0\|$, and $M$, and  the terms $K_1$ and $K_2$ are defined in (\ref{2.35}). From the Cauchy inequality, the Sobolev inequality, the Young inequality, and Corollary 3.1, we have
\begin{eqnarray}\label{3.30}
K_{1}&\leq&\eta\int_{0}^{t}\|u_{xx}(\tau)\|^{2}d\tau+C(\eta)\varepsilon^{4}\int_{0}^{t}\|u_{x}(\tau)\|^{\frac{2}{3}}\|v_{xx}(\tau)\|^{\frac{8}{3}}d\tau\nonumber\\
&\leq&\eta\int_{0}^{t}\|u_{xx}(\tau)\|^{2}d\tau+C(\eta)\int_{0}^{t}\|u_{x}(\tau)\|^{2}d\tau+C(\eta)\int_{0}^{t}\varepsilon^{4}\|v_{xx}(\tau)\|^{4}d\tau,
\end{eqnarray}
\begin{eqnarray}\label{3.31}
K_{2}&\leq&\eta\int_{0}^{t}\|u_{xx}(\tau)\|^{2}d\tau+C(\eta)\int_{0}^{t} \left(\|v_{x}(\tau)\|^{2}+\|u_{x}(\tau)\|\|u_{xx}(\tau)\|\|v_{x}(\tau)\|^{2}\right.\nonumber\\
&&+\left.\varepsilon^{4}\|v_{x}(\tau)\|\|v_{xx}(\tau)\|^{3}+\varepsilon^{4}\|v_{x}(\tau)\|^{4} \|v_{xx}(\tau)\|^{2}\right)d\tau\nonumber\\
&\leq&2\eta\int_{0}^{t}\|u_{xx}(\tau)\|^{2}d\tau+C(\eta)\int_{0}^{t}\left(\|v_{x}(\tau)\|^{2}+\sup_{0\leq \tau\leq t}\left\{{\|v_{x}(\tau)\|^{4}}\right\}\|u_{x}(\tau)\|^{2}\right.\nonumber\\
&&+\left.\varepsilon^{2}\|v_{x}(\tau)\|_{1}^{2}\cdot\varepsilon^{2}\|v_{xx}(\tau)\|^{2}+\varepsilon^{2}\|v_{xx}(\tau)\|^{2}\sup_{0\leq \tau\leq t}\{\varepsilon^{2}\|v_{x}(\tau)\|^{4}\}\right)d\tau\\
&\leq&2\eta\int_{0}^{t}\|u_{xx}(\tau)\|^{2}d\tau+C(\eta)\int_{0}^{t}\left[\|(v _{x}(\tau), u_x(\tau))\|^{2}+(\varepsilon^{2}\|v_{x}(\tau)\|_{1}^{2}+1)\varepsilon^{2}\|v_{xx}(\tau)\|^{2}\right]d\tau.\nonumber
\end{eqnarray}
Then inserting  (\ref{3.30})-(\ref{3.31}) into (\ref{2.29}), and using the Gronwall inequality and Corollary 3.1, we can obtain (\ref{3.28}). This finishes the proof of Lemma 3.3.

 For the estimate on $\|v_{xx}(t)\|$, we have
\begin{Lemma}
If the conditions   of  Proposition 3.2 hold, then  there exists a positive constant
$C_{31}(M)$ depending only on $\nu,\gamma,\underline{V}, \overline{V}$, $\|v_{0}-1\|_2$, $\|u_0\|_1$, and $M$, such that for $t\in[0,T]$,
\begin{eqnarray}\label{3.32}
\|v_{xx}(t)\|^{2}+\int_{0}^{t}\|v_{xx}(\tau)\|^{2}d\tau+\varepsilon^{2}\int_{0}^{t}\|v_{xxx}(\tau)\|^{2}d\tau\leq C_{31}(M)\left(\|v_{0}-1\|_{2}^{2}+\|u_{0}\|_{1}^{2}\right).
\end{eqnarray}
\end{Lemma}
\noindent{\bf Proof.}~~Using the same method as the proof of (\ref{2.40}), we have
\begin{eqnarray}\label{3.33}
&&\|v_{xx}(t)\|^{2}+\int_{0}^{t}\|v_{xx}(\tau)\|^{2}d\tau+\varepsilon^{2}\int_{0}^{t}\|v_{xxx}(\tau)\|^{2}d\tau\nonumber\\
&\leq&C_{32}(M)\left(\|(v_{0xx},u_{0x})\|^{2}+\int_{0}^{t}\|u_{xx}(\tau)\|^{2}d\tau+\|u_x(t)\|^2+\int_{0}^{t}\int_{\mathbb{R}}(|K_{3}|+|K_{4}|)dxd\tau\right),
\end{eqnarray}
where $C_{32}(M)$ is a positive constant depending only on $\nu,\gamma,\underline{V}, \overline{V}$, $\|v_{0}-1\|_1$, $\|u_0\|$, and $M$, and  the terms $K_3$ and $K_4$ are defined in (\ref{2.39}).
Similar to (\ref{3.30})-(\ref{3.31}), we have
\begin{eqnarray}\label{3.34}
C_{32}(M)\int_{0}^{t}\int_{\mathbb{R}}|K_{3}|dxd\tau&\leq&C(M)\int_{0}^{t}\int_{\mathbb{R}}
\left(|v_{x}^{2}v_{xx}|+|v_{xx}^{2}u_{x}|+|v_{x}v_{xx}u_{xx}|+|v_{x}^{2}v_{xx}u_{x}|\right)dxd\tau\nonumber\\
&\leq&\eta\int_{0}^{t}\|v_{xx}(\tau)\|^{2}d\tau+C(\eta,M)\int_{0}^{t}\|v_{x}(\tau)\|^{3}\|v_{xx}(\tau)\|d\tau\nonumber\\
&&+C(\eta,M)\int_{0}^{t}\left(\|u_{x}(\tau)\|_1^2+\|u_{x}(\tau)\|^{2}\|v_{x}(\tau)\|^{2}\right)\|v_{xx}(\tau)\|^{2}d\tau\nonumber\\
&&+C(M)\int_{0}^{t}\|v_{x}(\tau)\|^{\frac{1}{2}}\|v_{xx}(\tau)\|^{\frac{3}{2}}\|u_{xx}(\tau)\|d\tau\\
&\leq&2\eta\int_{0}^{t}\|v_{xx}(\tau)\|^{2}d\tau+C(\eta,M)\int_{0}^{t}(\|u_{x}(\tau)\|_1^{2}\|v_{xx}(\tau)\|^{2}+\|v_x(\tau)\|^{2})d\tau,\nonumber
\end{eqnarray}
\begin{eqnarray}\label{3.35}
&&C_{32}(M)\int_{0}^{t}\int_{\mathbb{R}}|K_{4}|dxd\tau\nonumber\\
&\leq&\frac{\varepsilon^{2}}{2}\int_{0}^{t}\|v_{xxx}(\tau)\|^{2}d\tau
+C(M)\varepsilon^{2}\int_{0}^{t}\int_{\mathbb{R}}\left(|v_{x}^{2}v_{xx}^{2}|+|v_{x}^{6}|\right)dxd\tau\nonumber\\
&\leq&\frac{\varepsilon^{2}}{2}\int_{0}^{t}\|v_{xxx}(\tau)\|^{2}d\tau
+C(M)\varepsilon^{2}\int_{0}^{t}\left(\|v_{x}(\tau)\|_{L^\infty}^{2}\|v_{xx}(\tau)\|^{2}+\|v_{x}(\tau)\|_{L^\infty}^{4}\|v_{x}(\tau)\|^{2}\right)dxd\tau\nonumber\\
&\leq&\frac{\varepsilon^{2}}{2}\int_{0}^{t}\|v_{xxx}(\tau)\|^{2}d\tau
+C(M)\varepsilon^{2}\int_{0}^{t}\left(\|v_{x}(\tau)\|\|v_{xx}(\tau)\|^{3}+\sup_{0\leq\tau\leq t}\{\|v_{x}(\tau)\|^{4}\}\|v_{xx}(\tau)\|^{2}\right)d\tau\nonumber\\
&\leq&\frac{\varepsilon^{2}}{2}\int_{0}^{t}\|v_{xxx}(\tau)\|^{2}d\tau
+C(M)\left(\varepsilon^{2}\int_{0}^{t}\|v_{x}(\tau)\|_{1}^{2}d\tau+\varepsilon^{2}\int_{0}^{t}\|v_{xx}(\tau)\|^{4}d\tau\right).
\end{eqnarray}
 Then (\ref{3.32}) follows immediately by combining (\ref{3.33})-(\ref{3.35}), Corollary 3.1, Lemma 3.3, and the Gronwall inequality.  The proof of Lemma 3.4 is completed.

By using the same argument as above, we can also obtain
\begin{Lemma}
If the conditions  of  Proposition 3.2 hold,  then there exists a positive constant
$C_{33}(M)$ depending only on $\nu,\gamma,\underline{V}, \overline{V}$, $\|v_{0}-1\|_5$, $\|u_0\|_4$, and $M$, such that for $t\in[0,T]$,
\begin{eqnarray}\label{3.36}
&&\quad\|v_{xxx}(t)\|_2^{2}+\|u_{xx}(t)\|_2^{2}+\varepsilon^2\|v_{xxx}(t)\|_2^{2}+\int_{0}^{t}\left(\|(v_{xxx},u_{xxx})(\tau)\|_2^{2}+\varepsilon^2\|v_{xxxx}(\tau)\|_2^{2}\right)d\tau\nonumber\\
&&\leq C_{33}(M)\left(\|v_{0}-1\|_{5}^{2}+\|u_{0}\|_{4}^{2}\right).
\end{eqnarray}
\end{Lemma}

\noindent{\bf Proof of Proposition 3.2.}~~ Proposition 3.2 follows immediately from Lemmas 3.1-3.5, and Corollary 3.1.
\subsection{Proof of Theorem 1.2}
Based on Propositions 3.1-3.2, we now use a  continuation argument to extend  the local  solution of  the Cauchy problem (\ref{1.5})-(\ref{1.6}) to be a global one, that is $T=+\infty$. And as a by-product,  we can obtain the uniform-in-time positive lower and upper bounds of the specific volume $v(t,x)$. Precisely, the proof of Theorem 1.2 can be  divided into the following four steps.\\
\noindent{\bf Step 1.}  First, for notational simplicity, we denote $\mathcal{E}(t):=\|v(t)-1\|^2_{5}+\|u(t)\|^2_{4}+\varepsilon^2\|v_x(t)\|^2_{4}$ and  $\mathcal{D}(t)=\|(v_x,u_x)(t)\|_4^2+\varepsilon^2\|v_{xx}(t)\|^2_{4}$.  Let $T_1=\frac{32}{\gamma}(\gamma+3)^2C_{24}^2C_{35}^2C_{27}$, where the constants $C_{24}$, $C_{27}$ and $C_{35}$ are  defined in (\ref{3.5}), (\ref{3.6}) and (\ref{3.47}), respectively. Then by the local existence result Proposition 3.1, the Cauchy problem (\ref{1.5})-(\ref{1.6}) has a unique solution  $(v(t,x),u(t,x))\in X_4(0,\tilde{t}_1;\frac{\underline{V}}{2},2\overline{V})$ with $\tilde{t}_1=\min\{T_1, \tilde{t}_0(\underline{V}, \tilde{N}_0)\}$,  which satisfies
 \begin{equation}\label{3.37}
  \sup_{[0,\tilde{t}_1]}\mathcal{E}(t)+\int_{0}^{\tilde{t}_1}\mathcal{D}(\tau)d\tau \leq \tilde{b}\tilde{N}_0^2,
\end{equation}
where  $\tilde{N}_0=(\|v_0-1\|^2_{5}+\|u_0\|^2_{4})^{\frac{1}{2}}$, and $\tilde{b}>1$ is a constant depending only on $\underline{V},\overline{V}$.

Then there exists a sufficiently small constant $\varepsilon_1>0$  such that
$\varepsilon_1^2(2\overline{V})^{\gamma-1}\tilde{b}^2\tilde{N}_0^4<1$. Consequently,  if $0<\varepsilon\leq\varepsilon_1$, we have from Proposition 3.2 that for all $(t,x)\in[0,\tilde{t}_1]\times\mathbb{R}$,
 \begin{equation}\label{3.38}
C_{24}^{-1}\leq v(t,x)\leq  \left(\sup_{\mathbb{R}}v_0(x)+C_{25}\right)e^{C(\nu,\gamma)\tilde{t}_1}:=v_M(\tilde{t}_1)\leq v_M(T_1),
\end{equation}
\begin{equation}\label{3.39}
  \mathcal{E}(t)+\int_{0}^{t}\mathcal{D}(\tau)d\tau \leq C_{26}(v_M(T_1))\tilde{N}_0^2.
\end{equation}

Next, taking $(v(\tilde{t}_1,x),u(\tilde{t}_1,x))$ as initial data,  then by Proposition 3.1 again, we can extend the local solution $(v(t,x),u(t,x))$ to the time interval $[0, \tilde{t}_1+\tilde{t}_2]$ with
\[\tilde{t}_2=\min\left\{T_1-\tilde{t}_1, \,\,\tilde{t}_0\left(C_{24}^{-1}, \sqrt{C_{26}(v_M(T_1))\tilde{N}_0^2}\right)\right\},\]
and $(v(t,x),u(t,x))\in X_4\left(0,\tilde{t}_1+\tilde{t}_2;\frac{C_{24}^{-1}}{2}, 2v_M(T_1)\right)$, which satisfies that
 \begin{equation}\label{3.40}
   \sup_{[0,\tilde{t}_1+\tilde{t}_2]}\mathcal{E}(t)+\int_{0}^{\tilde{t}_1+\tilde{t}_2}\mathcal{D}(\tau)d\tau \leq2\tilde{b}_1C_{26}(v_M(T_1))\tilde{N}_0^2
\end{equation}
with $\tilde{b}_1=\tilde{b}(C_{24}^{-1}, v_M(T_1))>1$ being  a constant depending only on $C_{24}^{-1}$ and $v_M(T_1)$.

Then there exists a sufficiently small constant $\varepsilon_2>0$  such that $$\varepsilon_2^2(2v_M(T_1))^{\gamma-1}(2\tilde{b}_1C_{26}(v_M(T_1))\tilde{N}_0^2)^2<1,$$
thus if $0<\varepsilon\leq\min\{\varepsilon_1,\varepsilon_2\}$, we have from Proposition 3.2 that for all $(t,x)\in[0,\tilde{t}_1+\tilde{t}_2]\times\mathbb{R}$,
 \begin{equation}\label{3.41}
C_{24}^{-1}\leq v(t,x)\leq  \left(\sup_{\mathbb{R}}v_0(x)+C_{25}\right)e^{C(\nu,\gamma)(\tilde{t}_1+\tilde{t}_2)}:=v_M(\tilde{t}_1+\tilde{t}_2)\leq v_M(T_1),
\end{equation}
\begin{equation}\label{3.42}
  \mathcal{E}(t)+\int_{0}^{t}\mathcal{D}(\tau)d\tau \leq C_{26}(v_M(T_1))\tilde{N}_0^2.
\end{equation}

Now taking $(v(\tilde{t}_1+\tilde{t}_2,x), u(\tilde{t}_1+\tilde{t}_2,x))$ as initial data and  employing  Proposition 3.1 again, if $0<\varepsilon\leq\min\{\varepsilon_1,\varepsilon_2\}$, then the local solution $(v(t,x), u(t,x))$  can be extended to the time interval  $[0,\tilde{t}_1+2\tilde{t}_2]$. By repeating the above procedure, then the local solution $(v(t,x), u(t,x))$  is extended to the time interval $[0,T_1]$. Moreover, the solution $(v(t,x), u(t,x))$ satisfies (\ref{3.41})-(\ref{3.42}) for all $(t,x)\in[0,T_1]\times\mathbb{R}$.

\noindent{\bf Step 2.}  We derive from Lemma 2.1 and (\ref{3.6}) that  for all $t\in[0,T_1]$,
\begin{equation}\label{3.43}
     \begin{split}
     \int_{\mathbb{R}}\left(\Phi(v)+\frac{u^2}{2}+\frac{\varepsilon^2v_x^2}{2v^4}\right)dx
     +\int_0^t\int_{\mathbb{R}}\frac{2\nu u_x^2}{v^{2}}dxd\tau
\leq C(\underline{V}, \bar{V})\left(\|v_0-1\|_1^2+\|u_0\|^2\right):=C_{34},
     \end{split}
\end{equation}
\begin{equation}\label{3.44}
    \int_{\mathbb{R}}\dfrac{\nu v_x^2}{v^{4}}dx+\int_0^t\int_{\mathbb{R}}\dfrac{\gamma v_x^2}{v^{\gamma+3}}dxd\tau\leq C_{27},
\end{equation}
where $C_{27}$ is a positive constant defined in (\ref{3.6}).
 From (\ref{3.44}), we have
\begin{equation}\label{3.45}
\int^{T_1}_{\frac{T_1}{2}}\left\|\frac{\sqrt{\gamma}v_x}{v^{s}}(\tau)\right\|^2d\tau\leq C_{27}
\end{equation}
with $s=\frac{\gamma+3}{2}$. Then from the definition of $T_1$ and (\ref{3.45}), we can  conclude  that  there exists a number $t_0^\prime\in[\frac{T_1}{2}, T_1]$ such that
\begin{equation}\label{3.46}
\left\|\frac{\sqrt{\gamma}v_x}{v^{s}}(t_0^\prime)\right\|\leq \frac{\sqrt{\gamma}}{4(\gamma+3)C_{24}C_{35}}.
\end{equation}

Since
 \begin{displaymath}
   \lim_{v\rightarrow C_{24}^{-1}}\frac{ \left(\frac{1}{v}-1\right)^2}{\Phi(v)}=c_1, \quad  \lim_{v\rightarrow +\infty}\frac{ \left(\frac{1}{v}-1\right)^2}{\Phi(v)}=0,
  \end{displaymath}
 where the function $\Phi(v)$ is defined in Lemma 2.1, and $c_1$ is a positive constant depending only on $C_{24}^{-1}$,  it holds that
 $\left(\frac{1}{v}-1\right)^2\leq c_2 \Phi(v)$
  for any $C_{24}^{-1}\leq v\leq+\infty$,  here $c_2$ is a positive constant depending only on $C_{24}^{-1}$. This together with (\ref{3.43}),  and (\ref{3.41}) with $(t,x)\in[0,T_1]\times\mathbb{R}$   implies that
  \begin{eqnarray}\label{3.47}
      \int_{\mathbb{R}}\left| v^{-s}-1\right|^2dx &\leq& \int_{\mathbb{R}}\left|\frac{1}{v}-1\right|^2s^2\left(\theta\frac{1}{v}+(1-\theta)\right)^{2(s-1)}dx\nonumber\\
      &\leq&  2^{2(s-1)}s^2\left(\left\|\frac{1}{v}\right\|_{L^\infty_{T_1,x}}^{2(s-1)}+1\right)c_2 \int_{\mathbb{R}}\Phi(v(t,x))dx \nonumber\\
         &\leq& 2^{2(s-1)}s^2((C_{24})^{2(s-1)}+1)C_{34}:= C^2_{35},
  \end{eqnarray}
where $\theta\in[0,1]$ is a constant.

On the other hand, we can deduce that
\begin{eqnarray*}
      \left(v^{-s}-1\right)^2 &=& \int_{-\infty}^x 2\left(v^{-s}-1\right)\left(v^{-s}-1\right)_ydy\nonumber\\
      &\leq& 2s\left|\int_{\mathbb{R}}\left(v^{-s}-1\right)\frac{v_x}{v^{s+1}}dx\right|\nonumber\\
      &\leq&2s \left\|\frac{1}{v}\right\|_{L^\infty_{T_1,x}}\|(v^{-s}-1)(t)\|\left\|\frac{v_x}{v^{s}}(t)\right\|\\
      &\leq&2sC_{24}C_{35}\left\|\frac{v_x}{v^{s}}(t)\right\|,\quad \forall\, (t,x)\in[0,T_1]\times\mathbb{R},
 \end{eqnarray*}
which together with (\ref{3.46}) gives rise to
\begin{eqnarray}\label{3.48}
      \left|v^{-s}(t_0^\prime,x)-1\right|\leq\sqrt{\frac{2s}{\sqrt{\gamma}}C_{24}C_{35}}\left\|\frac{\sqrt{\gamma}v_x}{v^{s}}(t_0^\prime)\right\|^{\frac{1}{2}}\leq\frac{1}{2}.
 \end{eqnarray}
Therefore, we have from (\ref{3.41}) and (\ref{3.48}) that
\begin{eqnarray}\label{3.49}
      C_{24}^{-1}\leq v(t_0^\prime,x)\leq2^{\frac{1}{s}}=2^{\frac{2}{\gamma+3}},\quad \forall \,x \in\mathbb{R}.
 \end{eqnarray}

\noindent{\bf Step 3.} Now, we take $(v(t^\prime_0,x),u(t^\prime_0,x))$ as initial data,  then we have from Proposition 3.1 that  the local solution $(v(t,x),u(t,x))$ can be extended to the time interval $[0, t^\prime_0+t^\prime_1]$ with
\[t^\prime_1=\min\left\{T_1, \,\,\tilde{t}_0\left(C_{24}^{-1}, \sqrt{C_{26}(v_M(T_1))\tilde{N}_0^2}\right)\right\},\]
and $(v(t,x),u(t,x))\in X_4\left(0,t^\prime_0+t^\prime_1;\frac{C_{24}^{-1}}{2}, 2v_M(T_1)\right)$, which satisfies that
 \begin{equation}\label{3.50}
  \sup_{[0,t^\prime_0+t^\prime_1]}\mathcal{E}(t)+\int_{0}^{t^\prime_0+t^\prime_1}\mathcal{D}(\tau)d\tau \leq2\tilde{b}_2C_{26}(v_M(T_1))\tilde{N}_0^2
\end{equation}
with $\tilde{b}_2=\tilde{b}(C_{24}^{-1}, v_M(T_1))>1$ being  a constant depending only on $C_{24}^{-1}$ and $v_M(T_1)$.

Then there exists a sufficiently small constant $\varepsilon_3>0$  such that
 \[
  \varepsilon_3^2(2v_M(T_1))^{\gamma-1}(2\tilde{b}_2C_{26}(v_M(T_1))\tilde{N}_0^2)^2<1.
\]
Consequently, if $0<\varepsilon\leq\min\{\varepsilon_1,\varepsilon_2,\varepsilon_3\}$, we have from Proposition 3.2 that for all $(t,x)\in[0,t^\prime_0+t^\prime_1]\times\mathbb{R}$,
 \begin{eqnarray}\label{3.51}
C_{24}^{-1}\leq v(t,x)&\leq& \left(\sup_{\mathbb{R}}v(t^\prime_0,x)+C_{25}\right)e^{C(\nu,\gamma)t^\prime_1}\nonumber\\
&\leq&  \left(2^{\frac{2}{\gamma+3}}+C_{25}\right)e^{C(\nu,\gamma)t^\prime_1}:=\hat{v}_M(t^\prime_1)\leq \hat{v}_M(T_1),
\end{eqnarray}
\begin{equation}\label{3.52}
\mathcal{E}(t)+\int_{0}^{t}\mathcal{D}(\tau)d\tau \leq C_{26}(\hat{v}_M(T_1))\tilde{N}_0^2.
\end{equation}

Next, taking $(v(t^\prime_0+t^\prime_1,x),u(t^\prime_0+t^\prime_1,x))$ as initial data and using  Proposition 3.1 again, then  the local solution $(v(t,x),u(t,x))$ can be extended to the time interval $[0, t^\prime_0+t^\prime_1+t^\prime_2]$ with
\[t^\prime_2=\min\left\{T_1-t^\prime_1, \,\,\tilde{t}_0\left(C_{24}^{-1}, \sqrt{C_{26}(\hat{v}_M(T_1))\tilde{N}_0^2}\right)\right\},\]
and $(v(t,x),u(t,x))\in X_4\left(0,t^\prime_0+t^\prime_1+t^\prime_2;\frac{C_{24}^{-1}}{2}, 2\hat{v}_M(T_1)\right)$, which satisfies that
 \begin{eqnarray}\label{3.53}
  \sup_{[0,  t^\prime_0+t^\prime_1+t^\prime_2]}\mathcal{E}(t)+\int_{0}^{ t^\prime_0+ t^\prime_1+t^\prime_2}\mathcal{D}(\tau)d\tau
  \leq 2\tilde{b}_3 C_{26}(\hat{v}_M(T_1))\tilde{N}_0^2.
\end{eqnarray}
Here $\tilde{b}_3=\tilde{b}(C_{24}^{-1}, \hat{v}_M(T_1))>1$ being  a constant depending only on  $C_{24}^{-1}$ and $\hat{v}_M(T_1)$.

Then there exists a sufficiently small constant $\varepsilon_4>0$  such that
 \[
  \varepsilon_4^2(2\hat{v}_M(T_1))^{\gamma-1}(2\tilde{b}_3C_{26}(\hat{v}_M(T_1))\tilde{N}_0^2)^2<1.
\]
Consequently, if $0<\varepsilon\leq\min\{\varepsilon_1,\varepsilon_2,\varepsilon_3,\varepsilon_4\}$, it follows from Proposition 3.2 that for all $(t,x)\in[0,t^\prime_0+t^\prime_1+t^\prime_2]\times\mathbb{R}$,
 \begin{equation}\label{3.54}
C_{24}^{-1}\leq v(t,x)\leq  \left(\sup_{\mathbb{R}}v(t^\prime_0,x)+C_{25}\right)e^{C(\alpha,\gamma)(t^\prime_1+t^\prime_2)}\leq\hat{v}_M(t^\prime_1+t^\prime_2)\leq \hat{v}_M(T_1),
\end{equation}
\begin{equation}\label{3.55}
\mathcal{E}(t)+\int_{0}^{t}\mathcal{D}(\tau)d\tau
\leq C_{26}(\hat{v}_M(T_1))\tilde{N}_0^2.
\end{equation}

Now taking  $(v(t^\prime_0+t^\prime_1+t^\prime_2,x),u(t^\prime_0+t^\prime_1+t^\prime_2,x))$ as initial data, then by Proposition 3.1 again, the local solution  can be extended to the time interval  $[0,t^\prime_0+t^\prime_1+2t^\prime_2]$ provided that $0<\varepsilon\leq\min\{\varepsilon_1,\varepsilon_2,\varepsilon_3,\varepsilon_4\}$. Repeating the above argument, the local solution $(v(t,x), u(t,x))$   can be extended by finite steps to the time interval $[0,t^\prime_0+T_1]$. Moreover, the solution $(v(t,x), u(t,x))$ satisfies (\ref{3.54})-(\ref{3.55}) for all $(t,x)\in[0,t^\prime_0+T_1]\times\mathbb{R}$. Since $t^\prime_0+T_1\geq\frac{T_1}{2}+T_1=\frac{3T_1}{2}$, we have thus extended the local solution to the time interval $[0,\frac{3T_1}{2}]$.

\noindent{\bf Step 4.} Since (\ref{3.43})-(\ref{3.44}) also hold for all $t\in[0,\frac{3T_1}{2}]$, we can choose a time $t_0^{\prime\prime}\in[t_0^\prime+\frac{T_1}{2},t_0^\prime+T_1]$ such that (\ref{3.49}) holds with $t_0^{\prime}$ replaced by $t_0^{\prime\prime}$. Hence by repeating the same  argument as that in Step 3,  we can  extend the local solution to the time interval $[0,t_0^{\prime\prime}+T_1]$ provided that $0<\varepsilon\leq\varepsilon_0:=\min\{\varepsilon_1,\varepsilon_2,\varepsilon_3,\varepsilon_4\}$. Since $t^{\prime\prime}_0+T_1\geq t^{\prime}_0+\frac{T_1}{2}+T_1\geq2T_1$, we have  extended the local solution $(v(t,x), u(t,x))$ to the time interval $[0,2T_1]$, and $(v(t,x), u(t,x))$ satisfies (\ref{3.54})-(\ref{3.55}) for all $(t,x)\in[0,2T_1]\times\mathbb{R}$.

By repeating the above procedure, the local solution can be extended  to a global one provided that $0<\varepsilon\leq\varepsilon_0$. Moreover, the estimates (\ref{3.54})-(\ref{3.55}) hold for all $(t,x)\in[0,+\infty)\times\mathbb{R}$. Letting
\[C_2:=\max\{C_{24}, \hat{v}_M(T_1)\},\quad C_3:=C_{26}(\hat{v}_M(T_1)),\]
then (\ref{1.13}) and (\ref{1.14}) follow immediately from  (\ref{3.54}) and (\ref{3.55}) with $(t,x)\in[0,+\infty)\times\mathbb{R}$, respectively. The proof of the large-time behavior of solutions (\ref{1.15}) is similar to that of  (\ref{1.9}), and thus is omitted here for brevity. This finishes  the proof of Theorem 1.2.

\section{The vanishing dispersion limit}
\setcounter{equation}{0}
In this section, we shall present the proof of  Theorem 1.3, which is divided into the following two subsections.
\subsection{The  limit as the Plank constant $\varepsilon$ vanishes}
The aim of this subsection is to pass the limit   $\varepsilon\rightarrow0$ in the compressible quantum Navier-Stokes equations (\ref{1.5}). To this end, let $(v^{\varepsilon}(t,x), u^{\varepsilon}(t,x))$ be the unique classical solution obtained in Theorem 1.2. Then we have from (\ref{1.14}) that
\begin{eqnarray}\label{4.1}
&&\|v^\varepsilon(t)-1\|^2_{5}+\|u^\varepsilon(t)\|^2_{4}+\varepsilon^2\|v^\varepsilon_x(t)\|^2_{4}+
\int_0^t\left(\|(u^\varepsilon_x(s),v^\varepsilon_x(s))\|^2_{4}+\varepsilon^2\|v^\varepsilon_{xx}(s)\|^2_{4}\right)ds\nonumber\\
&&\leq C_3\left(\|v_0-1\|^2_{5}+\|u_0\|^2_{4}\right), \quad \forall\, t>0,
\end{eqnarray}
where $C_3$ is a positive constant depending only on $\nu, \gamma,\underline{V}, \overline{V},  \|v_0-1\|_{5}$, and $\|u_0\|_4$.

(\ref{4.1}) implies that there exists a subsequence of $(v^{\varepsilon}(t,x), u^{\varepsilon}(t,x))$, still denoted by $(v^{\varepsilon}(t,x), u^{\varepsilon}(t,x))$, such that as $\varepsilon\rightarrow0$,
\begin{eqnarray}\label{4.2}
&&v^{\varepsilon}(t,x)\rightharpoonup v^0(t,x),\quad \mbox{weakly-$\ast$  in}~ L^{\infty}(0,T;H^{5}(\mathbb{R})), \\[2mm]
&&u^{\varepsilon}(t,x)\rightharpoonup u^0(t,x), \quad \mbox{weakly-$\ast$ in}~L^{\infty}(0,T;H^{4}(\mathbb{R}))\label{4.3}
\end{eqnarray}
for any fixed $0<T\leq\infty$. On the other hand, (\ref{4.1}) together with $(\ref{1.5})$ yields
\begin{eqnarray*}
&&\|v_{t}^{\varepsilon}(t)\|_{3}=\|u_{x}^{\varepsilon}(t)\|_{3}\leq C_{3}\left(\|v_{0}-1\|_{5}^{2}+\|u_{0}\|_{4}^{2}\right),\\[2mm]
&&\|u_{t}^{\varepsilon}(t)\|_{2}\leq C_{3}\left(\|v_{0}-1\|_{5}^{2}+\|u_{0}\|_{4}^{2}\right).
\end{eqnarray*}
Therefore, $v_{t}^{\varepsilon}(t,x)\in L^{\infty}(0,T;H^{3}(\mathbb{R}))$ and $u_{t}^{\varepsilon}(t,x)\in L^{\infty}(0,T;H^{2}(\mathbb{R}))$.

Since $H^{5}{(\mathbb{R})}\hookrightarrow\hookrightarrow C_{b}^{4}(\mathbb{R})\subset H^{3}{(\mathbb{R})}$ and $H^{4}{(\mathbb{R})}\hookrightarrow\hookrightarrow C_{b}^{3}(\mathbb{R})\subset H^{2}{(\mathbb{R})}$, it follows from the Lions-Aubin lemma (see Lemma 5.2) that there exists a subsequence of $(v^{\varepsilon}(t,x), u^{\varepsilon}(t,x))$, still denoted by $(v^{\varepsilon}(t,x), u^{\varepsilon}(t,x))$, such that as $\varepsilon\rightarrow0$,
\begin{eqnarray}\label{4.4}
&&v^{\varepsilon}(t,x)\rightarrow v^{0}(t,x) ~\mbox{strongly in} ~C\left(0,T; H_{loc}^{5-s}(\mathbb{R})\cap C_{b}^{4}(\mathbb{R})\right), \\[2mm]
&&u^{\varepsilon}(t,x)\rightarrow u^{0}(t,x) ~\mbox{strongly in}~ C\left(0,T; H_{loc}^{4-s}(\mathbb{R})\cap C_{b}^{3}(\mathbb{R})\right)\label{4.5}
\end{eqnarray}
with $s\in(0,\frac{1}{2})$. By virtue of
\[\int_{0}^{t}\left(\|v_{x}^{\varepsilon}(\tau)\|_{4}^{2}+\varepsilon^{2}\|v_{xx}^{\varepsilon}(\tau)\|_{4}^{2}\right)d\tau\leq C_{3}\left(\|v_{0}-1\|_{3}^{2}+\|u_{0}\|_{3}^{2}\right),\,\,\forall\,t>0,\]
and $C_{2}^{-1}\leq v^{\varepsilon}(t,x)\leq C_{2}$ with $C_2$ being a positive constant depending only on $\nu, \gamma,\underline{V}, \overline{V},  \|v_0-1\|_{1}$, and $\|u_0\|$, we obtain
\begin{equation}\label{4.6}
\varepsilon^{2}\left(-\frac{v_{xx}^{\varepsilon}(t,x)}{\left(v^{\varepsilon}(t,x)\right)^{4}}
+\frac{2\left(v_{x}^{\varepsilon}(t,x)\right)^{2}}{\left(v^{\varepsilon}(t,x)\right)^{5}}\right)_{x}\rightarrow0, \quad \mbox{strongly in} ~ L^{2}(0,T;H_{loc}^{2}(\mathbb{R})).
\end{equation}
Moreover, from the continuity  equation $(\ref{1.5})_1$ and (\ref{4.1}), we have
\[\int_{0}^{T}\|v_{t}^{\varepsilon}(\tau)\|_{H_{loc}^{3}(\mathbb{R})}^{2}d\tau
=\int_{0}^{T}\|u_{x}^{\varepsilon}(\tau)\|_{H_{loc}^{3}(\mathbb{R})}^{2}d\tau\leq C_{3}\left(\|v_{0}-1\|_{5}^{2}+\|u_{0}\|_{4}^{2}\right),\]
which together with (\ref{4.2}) yields
\begin{eqnarray}\label{4.7}
v^{\varepsilon}_t(t,x)\rightharpoonup v_t^0(t,x),\quad \mbox{weakly  in}~ L^{2}(0,T;H^{3}_{loc}(\mathbb{R})), \quad \mbox{as}~ \varepsilon\rightarrow0.
\end{eqnarray}

Similarly, it also holds  that
\begin{eqnarray}\label{4.8}
u^{\varepsilon}_t(t,x)\rightharpoonup u_t^0(t,x),\quad \mbox{weakly  in}~ L^{2}(0,T;H^{2}_{loc}(\mathbb{R})),\quad \mbox{as}~ \varepsilon\rightarrow0.
\end{eqnarray}

By virtue of (\ref{4.4})-(\ref{4.8}), we see that the limit function $(v^0(t,x), u^0(t,x))$ satisfies the compressible Navier-Stokes equations  (\ref{1.16}) for almost $(x,t)\in[0,T]\times\mathbb{R}$. Passing the limit $\varepsilon\rightarrow0$ in the inequalities $C_{2}^{-1}\leq v^{\varepsilon}(t,x)\leq C_{2}$  and taking into account (\ref{4.4})-(\ref{4.5}), then  $(v^0(t,x), u^0(t,x))$ satisfies (\ref{1.18-1}) and the initial data (\ref{1.17}). Furthermore, in view of (\ref{4.1}) and the  lower semi-continuity properties of the \mbox{weak-$\ast$} and weak topology, $(v^0(t,x), u^0(t,x))$ satisfies (\ref{1.18}),  which implies that  $(v^0(t,x), u^0(t,x))$ is a classical solution  of the one-dimensional compressible Navier-Stokes equations (\ref{1.16}). The large time behavior (\ref{1.18-2}) is a direct consequence of (\ref{1.18}) and the system (\ref{1.16}), the proof of which is similar to that of (\ref{1.9}) and thus is omitted here.  This finishes the proof of the first part of Theorem 1.3, i.e., the limit $\varepsilon\rightarrow0$.

\subsection{Convergence rates}
Now we are ready to prove the convergence rate estimates (\ref{1.19})-(\ref{1.20}) in Theorem 1.3 by using the energy method. First, let $(v^\varepsilon(t,x), u^\varepsilon(t,x))$ and $(v^0(t,x), u^0(t,x))$  be the classical solutions of the Cauchy problems  (\ref{1.5})-(\ref{1.6}), and (\ref{1.16})-(\ref{1.17}), respectively. Define
\[V(t,x)=v^{\varepsilon}(t,x)-v^{0}(t,x), \quad U(t,x)=u^{\varepsilon}(t,x)-u^{0}(t,x),\]
then it follows from (\ref{1.5}) and (\ref{1.16}) that
\begin{eqnarray}\label{4.9}
\left\{\begin{array}{ll}
    V_{t}-U_{x}=0, \,\,t>0, x\in\mathbb{R},\\[2mm]
    U_{t}+\left[p(v^{\varepsilon})-p(v^{0})\right]_{x}
    =\displaystyle 2\nu\left(\frac{U_{x}}{(v^{\varepsilon})^{2}}\right)_{x}
    +F
\end{array}\right.
\end{eqnarray}
with the initial data
\begin{eqnarray}\label{4.10}
(V(t,x),U(t,x))|_{t=0}=(0,0),\,\,\,x\in\mathbb{R},
\end{eqnarray}
where the nonlinear function $F$ is given by
\[F=2\nu\left[u_{x}^{0}\left(\frac{1}{(v^{\varepsilon})^{2}}-\frac{1}{(v^{0})^{2}}\right)\right]_{x}+\varepsilon^{2}\left(-\frac{v_{xx}^{\varepsilon}}{(v^{\varepsilon})^{4}}+\frac{2(v_{x}^{\varepsilon})^{2}}{(v^{\varepsilon})^{5}}\right)_{x}.\]

We divide the analysis into the following three steps:

\noindent {\bf Step 1.}~~The estimate on $\|(V(t),U(t))\|$.

Multiplying $(\ref{4.9})_{1}$ by $p(v^{0})-p(v^{\varepsilon})$, $(\ref{4.9})_{2}$ by $U$, then adding them together and integrating the resultant equation over $\mathbb{R}$, we have
\begin{eqnarray}\label{4.11}
&&\frac{d}{dt}\int_{\mathbb{R}}\left[\Phi(v^{\varepsilon},v^{0})+\frac{U^{2}}{2}\right]dx+2\nu\int_{\mathbb{R}}\frac{U_{x}^{2}}{(v^{\varepsilon})^{2}}dx\nonumber\\
&=&-\int_{\mathbb{R}}v_{t}^{0}[p(v^{\varepsilon})-p(v^{0})-p'(v^{0})V]dx
+2\nu\int_{\mathbb{R}}\left[u_{x}^{0}\left(\frac{1}{(v^{\varepsilon})^{2}}-\frac{1}{(v^{0})^{2}}\right)\right]_{x}Udx\nonumber\\
&&+\varepsilon^{2}\int_{\mathbb{R}}\left(-\frac{v_{xx}^{\varepsilon}}{(v^{\varepsilon})^{4}}+\frac{2(v_{x}^{\varepsilon})^{2}}{(v^{\varepsilon})^{5}}\right)_{x}Udx\nonumber\\
&:=&Q_{1}+Q_{2}+Q_{3},
\end{eqnarray}
where the function $\Phi(v^{\varepsilon},v^{0})=p(v^{0})(v^{\varepsilon}-v^{0})-\int_{v^{0}}^{v^{\varepsilon}}p(s)ds $.

From the Cauchy inequality, the Sobolev inequality, the uniform-in-$\varepsilon$ estimates (\ref{4.1}), (\ref{1.18}), and the fact that $C^{-1}_2\leq v^\varepsilon(t,x), v^0(t,x)\leq C_2$, we have
\begin{eqnarray}\label{4.12}
Q_{1}&\leq& C\int_{\mathbb{R}}|u_{x}^{0}|\cdot|V^{2}|dx\leq C\|u_{x}^{0}\|_{L^{\infty}}\|V(t)\|^{2}\leq C\|V(t)\|^{2},\\[2mm]
Q_{2}+Q_{3}&\leq&\eta\|U_{x}(t)\|^{2}+C_{\eta}\int_{\mathbb{R}}\left((u_{x}^{0})^{2}V^{2}+\varepsilon^{4}\left(v_{xx}^{\varepsilon}\right)^{2}+\varepsilon^{4}\left(v_{x}^{\varepsilon}\right)^{4}\right)dx\nonumber\\
&\leq&\eta\|U_{x}(t)\|^{2}+C_{\eta}\left(\|u_{x}^{0}(t)\|_{L^{\infty}}^{2}\cdot\|V(t)\|^{2}+\varepsilon^{4}\|v^\varepsilon_{x}(t)\|_{L^{\infty}}^{2}\|(v_{x}^{\varepsilon},v_{xx}^{\varepsilon})(t)\|^{2}\right)\nonumber\\
&\leq&\eta\|U_{x}(t)\|^{2}+C_{\eta}\left(\|V(t)\|^{2}+\varepsilon^{4}\right).\label{4.13}
\end{eqnarray}

Substituting (\ref{4.12})-(\ref{4.13}) into (\ref{4.11}), and integrating the resultant equation with respect to $t$ over $[0,t]$, we obtain by the smallness of $\eta$ that
\begin{eqnarray*}
\|(V(t),U(t))\|^{2}+\int_{0}^{t}\|U_{x}(\tau)\|^{2}d\tau\leq C\int_{0}^{t}\|V(\tau)\|^{2}d\tau+C\varepsilon^{4}t,
\end{eqnarray*}
which together with the Gronwall inequality leads to
\begin{eqnarray}\label{4.14}
\|(V(t),U(t))\|^{2}\leq C_{36}\varepsilon^{4}te^{\lambda_{1}t},
\end{eqnarray}
where $C_{36}$ and $\lambda_1$ are positive constants depending only on $\nu, \gamma,\underline{V}, \overline{V},  \|v_0-1\|_{5}$, and $\|u_0\|_4$.

\noindent{\bf Step 2.}~~The estimate of $\|(V_x(t), U_x(t))\|$.

Taking $(\ref{4.9})_{1x}\times V_x+(\ref{4.9})_{2x}\times U_{x}$, and integrating  the resultant equation with respect to $x$ over $\mathbb{R}$, we have
\begin{eqnarray}\label{4.15}
&&\frac{1}{2}\frac{d}{dt}\int_{\mathbb{R}}\left(V_{x}^{2}+U_{x}^{2}\right)dx+2\nu\int_{\mathbb{R}}\frac{U_{xx}^{2}}{(v^{\varepsilon})^{2}}dx\nonumber\\
&=&\int_{\mathbb{R}}U_{xx}V_{x}dx-\int_{\mathbb{R}}U_{x}[p(v^{\varepsilon})-p(v^{0})]_{xx}dx
+2\nu\int_{\mathbb{R}}\left[u_{x}^{0}\left(\frac{1}{(v^{\varepsilon})^{2}}-\frac{1}{(v^{0})^{2}}\right)\right]_{xx}U_{x}dx\nonumber\\
&&+\varepsilon^{2}\int_{\mathbb{R}}\left(-\frac{v_{xx}^{\varepsilon}}{(v^{\varepsilon})^{4}}+\frac{2(v_{x}^{\varepsilon})^{2}}{(v^{\varepsilon})^{5}}\right)_{xx}U_{x}dx+4\nu\int_{\mathbb{R}}\frac{v_{x}^{\varepsilon}U_{x}U_{xx}}{(v^{\varepsilon})^{3}}dx\nonumber\\
&:=&Q_{4}+Q_{5}+Q_{6}+Q_{7}+Q_{8}.
\end{eqnarray}
Similar to the estimates of (\ref{4.12})-(\ref{4.13}), we have
\begin{eqnarray}\label{4.16}
\sum_{i=4}^8Q_i\leq\eta\|U_{xx}(t)\|^{2}+C_{\eta}\left(\|(V_x,U_{x})(t)\|^{2}+\|V(t)\|^{2}+\varepsilon^{4}\right).
\end{eqnarray}

Combining (\ref{4.15})-(\ref{4.16}) and integrating the resultant equation  with respect to $t$  over $[0,t]$, we derive from (\ref{4.14}), the smallness of $\eta$ and  the  Gronwall inequality that
\begin{eqnarray}\label{4.17}
\|(V_{x}(t),U_{x}(t))\|^{2}\leq C_{37}\varepsilon^{4}\left(1+t^{2}\right)e^{\lambda_{2}t},
\end{eqnarray}
where $C_{37}$ and $\lambda_2$ are positive constants depending only on $\nu, \gamma,\underline{V}, \overline{V},  \|v_0-1\|_{5}$, and $\|u_0\|_4$.

\noindent{\bf Step 3.}~~The estimate of $\|(V_{xx}(t), U_{xx}(t))\|_1$.

Similar to the Step 2, by performing $(\ref{4.9})_{1xx}\times V_{xx}$+$(\ref{4.9})_{2xx}\times U_{xx}$, and integrating  the resultant equation with respect to $x$ over $\mathbb{R}$, we obtain
\begin{eqnarray}\label{4.18}
&&\frac{1}{2}\frac{d}{dt}\int_{\mathbb{R}}\left(V_{xx}^{2}+U_{xx}^{2}\right)dx+\nu\int_{\mathbb{R}}\frac{U_{xxx}^{2}}{(v^{\varepsilon})^{2}}dx\nonumber\\
&=&\int_{\mathbb{R}}U_{xxx}V_{xx}dx-\int_{\mathbb{R}}[p(v^{\varepsilon})-p(v^{0})]_{xxx}U_{xx}dx
+4\nu\int_{\mathbb{R}}\left[2U_{xx}\frac{v_{x}^{\varepsilon}}{(v^{\varepsilon})^{3}}+U_{x}\left(\frac{v_{xx}^{\varepsilon}}{(v^{\varepsilon})^{3}}-3\frac{(v_{x}^{\varepsilon})^{2}}{(v^{\varepsilon})^{4}}\right)\right]U_{xx}dx\nonumber\\
&&+2\nu\int_{\mathbb{R}}\left[u_{x}^{0}\left(\frac{1}{(v^{\varepsilon})^{2}}-\frac{1}{(v^{0})^{2}}\right)\right]_{xxx}U_{xx}dx
+\varepsilon^{2}\int_{\mathbb{R}}\left(-\frac{v_{xx}^{\varepsilon}}{(v^{\varepsilon})^{4}}+\frac{2(v_{x}^{\varepsilon})^{2}}{(v^{\varepsilon})^{5}}\right)_{xxx}U_{xx}dx\nonumber\\
&:=&Q_{9}+Q_{10}+Q_{11}+Q_{12}+Q_{13}.
\end{eqnarray}
Similar to (\ref{4.12})-(\ref{4.13}), the terms $Q_i (i=9,10,11,12,13)$ can be controlled as follows.
\begin{eqnarray}\label{4.19}
&&\displaystyle\sum_{i=9}^{11}Q_i\leq\eta\|U_{xxx}(t)\|^{2}+C_{\eta}\left(\|V(t)\|_2^{2}+\|(U_{x},U_{xx})(t)\|^{2}\right),\\[2mm]
&&Q_{12}\leq\eta\|U_{xxx}(t)\|^{2}+C_{\eta}\|V(t)\|_{2}^{2},\\[3mm]\label{4.20}
&&Q_{13}\leq\eta\|U_{xxx}(t)\|^{2}+C_{\eta}\varepsilon^{4}\left(\|v_{xxxx}^{\varepsilon}(t)\|^{2}+\|v_{x}^{\varepsilon}(t)\|_{2}^{2}\right)\nonumber\\
&&\qquad\leq\eta\|U_{xxx}(t)\|^{2}+C_{\eta}\varepsilon^{4}.\label{4.21}
\end{eqnarray}

Inserting (\ref{4.19})-(\ref{4.21}) into (\ref{4.18}), and integrating the resultant equation with respect to $t$ over $[0,t]$,  then we have from (\ref{4.14}), (\ref{4.17}), the smallness of $\eta$ and the Gronwall inequality that
\begin{eqnarray}\label{4.22}
\|(V_{xx}(t),U_{xx}(t))\|^{2}\leq C_{38}\varepsilon^{4}(1+t^3)e^{\lambda_{3}t},
\end{eqnarray}
where $C_{38}$ and $\lambda_3$ are positive constants depending only on $\nu, \gamma,\underline{V}, \overline{V},  \|v_0-1\|_{5}$, and $\|u_0\|_4$.

By repeating the same argument as above, we also have
\begin{eqnarray}\label{4.23}
\|(V_{xxx}(t),U_{xxx}(t))\|^{2}\leq C_{39}\varepsilon^{4}(1+t^4)e^{\lambda_{4}t},
\end{eqnarray}
where $C_{39}$ and $\lambda_4$ are positive constants depending only on $\nu, \gamma,\underline{V}, \overline{V},  \|v_0-1\|_{5}$, and $\|u_0\|_4$.

Then (\ref{1.19})-(\ref{1.20}) follows from (\ref{4.14}), (\ref{4.17}), (\ref{4.22}) and (\ref{4.23}) immediately. This ends the proof of Theorem 1.3.

\section{Appendix}
\setcounter{equation}{0}

\subsection{Some useful lemmas}
To derive the estimate of $\displaystyle\int_{\mathbb{R}}\frac{v_x^2}{v^4}dx$, we shall use the  following  coercivity inequality which is  due to Germain and LeFloch \cite{Germain-LeFloch-2012}.
\begin{Lemma} [Theorem 2.1, \cite{Germain-LeFloch-2012}] Let $f(x): \mathbb{R}\rightarrow(0,+\infty)$ be a positive function which approaches a constant $\rho_{*}$  at $\pm\infty$, and $(f(x)-\rho_{*})\in H^2(\mathbb{R})$. Then the inequality
\begin{equation}\label{5.1}
\int_{\mathbb{R}}f^a(x)(f_{xx}(x))^2dx\geq\left(\frac{a-1}{3}\right)^2\int_{\mathbb{R}}f^{a-2}(x)(f_{x}(x))^4dx
\end{equation}
holds for any $a\neq1$, in which the constant in the right-hand side is optimal.
\end{Lemma}

The following Lions-Aubin lemma will be used in the vanishing  dispersion limit for the system (\ref{1.5}).
\begin{Lemma} [Lions-Aubin lemma, see \cite{Aubin,Simon-1990}] Assume that $X\subset E\subset Y$ are Banach spaces and $X\hookrightarrow\hookrightarrow E$. Then the following imbeddings are compact:
\begin{itemize}
\item[$(i)$]$\left\{\varphi\,\Big|\varphi\in L^q(0,T;X),\,\,\displaystyle\frac{\partial\varphi}{\partial t}\in L^1(0,T;Y)\right\}\hookrightarrow\hookrightarrow L^q(0,T;E)$, \,if $1\leq q\leq \infty$;
\item[$(ii)$]$\left\{\varphi\,\Big|\varphi\in L^\infty(0,T;X),\,\,\displaystyle\frac{\partial\varphi}{\partial t}\in L^r(0,T;Y)\right\}\hookrightarrow\hookrightarrow C([0,T]; E)$, \,if $1< r\leq \infty$.
\end{itemize}
\end{Lemma}

\subsection{The derivation of  system (\ref{2.22})}
In this subsection, we shall give the derivation of the system (\ref{2.22}), that is
\begin{eqnarray*}
\left\{\begin{array}{ll}
     v_{t}-\xi_{x}=\displaystyle c_{+}\left(\frac{v_{x}}{v^{2}}\right)_{x},\quad x\in\mathbb{R},\,\,t>0,\\[4mm]
      \xi_{t}+\displaystyle\frac{\gamma}{c_{+}}v^{1-\gamma}\xi=\displaystyle\frac{\gamma}{c_{+}}v^{1-\gamma}u+c_{-}\left(\frac{\xi_{x}}{v^{2}}\right)_{x}.
\end{array}\right.
\end{eqnarray*}

Indeed, let the effective velocity $\xi(t,x)=u(t,x)-c_{0}\displaystyle\frac{v_{x}(t,x)}{v^2(t,x)}$ with $c_0$ being a constant to be determined later. Then from  $(\ref{1.5})_{1}$, we have
\begin{equation}\label{5.2}
v_{t}-\xi_{x}=c_{0}\left(\frac{v_{x}}{v^{2}}\right)_{x}.
\end{equation}
On the other hand, we infer from $(\ref{1.5})_{2}$ that
\begin{eqnarray*}
\xi_{t}+p(v)_{x}&=&(2\nu-c_{0})\displaystyle\left(\frac{u_{x}}{v^{2}}\right)_x+\varepsilon^{2}\left(-\frac{v_{xx}}{v^4}+\frac{2v_{x}^{2}}{v^{5}}\right)_{x}\\
&=&(2\nu-c_{0})\left(\frac{\xi_{x}}{v^{2}}\right)_x+\left[\varepsilon^{2}-(2\nu-c_{0})c_{0}\right]\left(-\frac{v_{xx}}{v^4}+\frac{2v_{x}^{2}}{v^{5}}\right)_{x}.
\end{eqnarray*}

Let $\nu^2\geq\varepsilon^2$ and $\varepsilon^{2}-(2\nu-c_{0})c_{0}=0$, then we obtain
\begin{eqnarray*}
c_{0}=\nu+\sqrt{\nu^{2}-\varepsilon^2}, \quad \mbox{or}\quad  c_{0}=\nu-\sqrt{\nu^{2}-\varepsilon^2}.
\end{eqnarray*}

Setting $c_0=c_{+}:=\nu+\sqrt{\nu^{2}-\varepsilon^2}$ and $c_{-}:=\nu-\sqrt{\nu^{2}-\varepsilon^2}$, then the momentum equation  $(\ref{1.5})_{2}$ is transformed as
\begin{equation}\label{5.3}
\xi_{t}+p(v)_{x}=\displaystyle c_{-}\left(\frac{\xi_{x}}{v^{2}}\right)_{x}.
\end{equation}
Since
\begin{eqnarray*}
p(v)_{x}&=&-\gamma v^{-\gamma-1}v_{x}=-\gamma v^{-\gamma-1}\frac{v^{2}}{c_{+}}(u-\xi)\\
&=&-\frac{\gamma}{c_{+}}v^{1-\gamma}u+\frac{\gamma}{c_{+}}v^{1-\gamma}\xi,
\end{eqnarray*}
it follows from (\ref{5.3}) that
\begin{equation}\label{5.4}
\xi_{t}+\frac{\gamma}{c_{+}}v^{1-\gamma}\xi=\frac{\gamma}{c_{+}}v^{1-\gamma}u+c_{-}\left(\frac{\xi_{x}}{v^{2}}\right)_{x}.
\end{equation}
Then system (\ref{2.22}) is a direct consequence of (\ref{5.2}) and  (\ref{5.4}). This finishes the derivation of (\ref{2.22}).

\bigbreak

\begin{center}
{\bf Acknowledgement}
\end{center}
The research of Zhengzheng Chen is supported by the National Natural Science Foundation of China under contract 12171001 and the Support Program for Outstanding Young Talents in Universities of Anhui Province under contract gxyqZD2022007. The research of Huijiang Zhao is supported by the National Natural Science Foundation of China under contract 11731008. This work is also partially supported by
a grant from Science and Technology Department of Hubei Province under contract 2020DFH002.

\end{document}